\documentclass[preprint]{elsarticle}
\usepackage[T1]{fontenc}
\usepackage[utf8]{inputenc}
\usepackage[hyperfootnotes=false]{hyperref}
\usepackage{tabularx}
\usepackage{supertabular}
\usepackage{booktabs}
\usepackage{dsfont,amsmath}
\usepackage{color,graphics,graphicx}
\usepackage[table]{xcolor}
\usepackage{multirow}
\usepackage{enumitem}
\usepackage{lmodern}
\usepackage{setspace}
\usepackage{nicefrac}
\usepackage{amssymb,amsthm}
\usepackage{bbm}
\usepackage{a4wide}
\newcommand{\N}{\mathds{N}}
\newcommand{\KLEINO}{{\scriptstyle{\mathcal{O}}}}
\newcommand{\E}{\mathbb{E}}
\newcommand{\R}{\mathbb{R}}
\DeclareMathAccent{\verywidehat}{\mathord}{largesymbols}{'144}
\newcommand{\1}{\mathbf{1}}
\newcommand{\cov}{\mathbb{C}\textnormal{o\hspace*{0.02cm}v}}
\renewcommand{\P}{\mathbb{P}}
\newcommand{\var}{\mathbb{V}\hspace*{-0.05cm}\textnormal{a\hspace*{0.02cm}r}}
\newtheorem{theorem}{Theorem}
\newtheorem{proposition}[theorem]{Proposition}%
\newtheorem{lem}[theorem]{Lemma} 

\newtheorem{remark}{Remark}%
\newtheorem{assumption}{Assumption}
\allowdisplaybreaks[3]
\begin{document}
\renewcommand*{\thefootnote}{\fnsymbol{footnote}}

\title{Inference on the intraday spot volatility from high-frequency order prices with irregular microstructure noise}
\author[1]{Markus Bibinger\footnote{Financial support from the Deutsche Forschungsgemeinschaft (DFG) under grant 403176476 is gratefully acknowledged.}}

\address[1]{Faculty of Mathematics and Computer Science, Institute of Mathematics, Julius-Maximilians-Universit\"at W\"urzburg, markus.bibinger@uni-wuerzburg.de}

\begin{frontmatter}

\begin{abstract} We consider estimation of the spot volatility in a stochastic boundary model with one-sided microstructure noise for high-frequency limit order prices. Based on discrete, noisy observations of an It\^{o} semimartingale with jumps and general stochastic volatility, we present a simple and explicit estimator using local order statistics. We establish consistency and stable central limit theorems as asymptotic properties. The asymptotic analysis builds upon an expansion of tail probabilities for the order statistics based on a generalized arcsine law. In order to use the involved distribution of local order statistics for a bias correction, an efficient numerical algorithm is developed. We demonstrate the finite-sample performance of the estimation in a Monte Carlo simulation.
\end{abstract}

\begin{keyword}{arcsine law\sep limit order book\sep market microstructure\sep nonparametric boundary model\sep volatility estimation}\\[.25cm]
\end{keyword}

\end{frontmatter}
\thispagestyle{plain}
\onehalfspacing
\section{Introduction\label{sec:1}}
Time series of intraday prices are typically described as a discretized path of a continuous-time stochastic process. To have arbitrage-free markets the log-price process should be a semimartingale. Risk estimation based on high-frequency data at highest available observation frequencies has to take microstructure frictions into account. Disentangling these market microstructure effects from dynamics of the long-run price evolution has led to observation models with additive noise, see, for instance, \cite{hansen}, \cite{aitjac14} and \cite{Linton}. The market microstructure noise, modelling among other effects the oscillation of traded prices between bid and ask order levels in an electronic market, is classically a centred (white) noise process with expectation equal to zero. These models can explain many stylized facts of high-frequency data. Having available full limit order books including data of submissions, cancellations and executions of bid and ask limit orders, however, it is not clear which time series to consider at all. While challenging the concept of one price process it raises the question if the information can be exploited more efficiently, in particular to improve risk quantification. The considered stochastic boundary model for limit order prices of an order book has been discussed by \cite{BJR}, \cite{Liu} and Chapter 1.8 of \cite{bishwal}. It preserves the concept of an underlying efficient, semimartingale log-price which determines long-run price dynamics and an additive, exogenous noise which models market-specific microstructure frictions. Its key idea is that ask order prices should (in most cases) lie above the unobservable efficient price and bid prices below the efficient price. This leads to observation errors which are irregular in the sense of having non-zero expectation and a distribution with a lower- or upper-bounded support. Considering without loss of generality a model for (best) ask order prices, we obtain lower-bounded observation errors and use local minima for the estimation. Modelling (best) bid prices instead would yield a model with upper-bounded observation errors and local maxima could be used for an analogous estimation. Both can be combined in practice.

It is known that the statistical and probabilistic properties of models with irregular noise are very different than for regular noise and require other methods, see, for instance, \cite{bound1}, \cite{bound3} and \cite{bound2}. Therefore, our estimation methods and asymptotic theory are quite different compared to the market microstructure literature, while we can still profit from some of the techniques used there. In \cite{BJR} an estimator for the quadratic variation of a continuous semimartingale, that is, the integrated volatility, was proposed with convergence rate $n^{-1/3}$, based on $n$ discrete observations with one-sided noise. Optimality of the rate was proved in the standard asymptotic minimax sense. A main insight was that this convergence rate is better than the optimal rate, $n^{-1/4}$, under regular market microstructure noise. 

A recent strand of literature proposes structural, parametric market microstructure noise models incorporating information based on observed order book quantities as volume or trade types, see \cite{chaker}, \cite{li}, \cite{clinet1} and \cite{clinet2}. Splitting the noise in a parametric function of such quantities and residual noise, a plug-in estimation of integrated volatility can yield as well faster convergence rates than in the classical model with uninformative noise. While this effect of improved volatility estimation appears to be a similarity to our work, our viewpoint on market microstructure is quite distinct. We focus on a model with one-sided instead of centred noise, but we neither impose a parametric assumption on the noise, nor do we include additional trading information. Such refinements of a one-sided noise model, as discussed in the mentioned works for the centred noise model, might be of interest for future research when microstructure effects of bid and ask quotes are better understood. This could potentially further improve volatility estimation.

Inference on the spot volatility is one of the most important topics in the financial literature, see, for instance, \cite{mancinispot}, \cite{hautsch} and the references therein. In this work, we address spot volatility estimation for the model from \cite{BJR}. Using local minima over blocks of shrinking lengths $h_n\propto n^{-2/3}\propto (nh_n)^{-2}$, the resulting distribution of local minima in \cite{BJR} became involved and infeasible, such that a central limit theorem for the integrated volatility estimator could not be obtained. Our spot volatility estimator is related to a localized version of the estimator from \cite{BJR}, combined with truncation methods to eliminate jumps of the semimartingale. For the asymptotic theory, however, we follow a different approach choosing blocks of lengths $h_n$, where $h_n n^{2/3}\to\infty$ slowly. This allows to establish stable central limit theorems with the best achievable rate, arbitrarily close to $n^{-1/6}$, in the important special case of a semimartingale volatility. We exploit this to construct pointwise asymptotic confidence intervals. 

Although the asymptotic theory relies on block lengths that are slightly unbalanced by smoothing out the impact of the noise distribution on the distribution of local minima asymptotically, our numerical study demonstrates that the confidence intervals work well in realistic scenarios with block lengths which optimize the estimation performance. Robustness to different noise specifications is an advantage that is naturally implied by our approach. Our estimator is surprisingly simple, it is a local average of squared differences of block-wise minima times a constant factor which comes from moments of the half-normal distribution of the minimum of a Brownian motion over the unit time interval. This estimator is consistent. However, the stable central limit theorem at fast convergence rate requires a subtle bias correction which incorporates a more precise approximation of the asymptotic distribution of local minima. For that purpose, our analysis is based on a generalization of the arcsine law which gives the distribution of the proportion of time over some interval that a Brownian motion is positive. In order to compute numerically the bias-correction function, we introduce an efficient algorithm. Reducing local minima over many random variables to iterated minima of two random variables in each step combined with a convolution step, it can be interpreted as a kind of dynamic programming approach. It turns out to be much more efficient compared to the natural approximation by a Monte Carlo simulation and is a crucial ingredient of our numerical application. Our convergence rate is much faster than the optimal rate, $n^{1/8}$, for spot volatility estimation under regular noise, see \cite{hoffmann} and \cite{hautsch}. The main contribution of this work is to develop the probabilistic foundation for the asymptotic analysis of the estimator and to establish the stable central limit theorems, asymptotic confidence and a numerically practicable method. 

The methods and proof techniques to deal with jumps are inspired by the truncation methods pioneered in \cite{mancini} and summarized in Chapter 13 of \cite{JP}. Overall, the strategy and restrictions on jump processes are to some extent similar, while several details under irregular noise using order statistics are rather different compared to settings without noise or with regular centred noise as in \cite{bibwinkel2}. 

We introduce and further discuss our model in Section \ref{sec:2}. Section \ref{sec:3} presents estimation methods and Section \ref{sec:4} asymptotic results. The numerical application is considered in Section \ref{sec:5} and a Monte Carlo simulation study illustrates an appealing finite-sample performance of the method. All proofs are given in Section \ref{sec:6}.

\section{Model with lower-bounded, one-sided noise and assumptions\label{sec:2}}
Consider an It\^o semimartingale   
\begin{eqnarray}X_t&=&X_0+\int_0^t a_s\,ds+\int_0^t\sigma_s\,dW_s + \int_0^t \int_{\mathbb{R}}\delta(s,z)\1_{\{|\delta(s,z)|\leq 1 \}}(\mu-\nu)(ds,dz) \nonumber \\ && +\int_0^t \int_{\mathbb{R}}\delta(s,z)\1_{\{|\delta(s,z)|> 1 \}} \mu(ds,dz)\,,\,t\ge 0\,,\label{sm} \end{eqnarray}
with a one-dimensional standard Brownian motion $(W_t)$, defined on some filtered probability space $(\Omega^X,\mathcal{F}^X,(\mathcal{F}^X_t),\mathbb{P}^X)$. For the drift process $(a_t)$, and the volatility process $(\sigma_t)$, we impose the following quite general assumptions.
\begin{assumption}\label{sigma}
The processes $(a_t)_{t\ge 0}$  and $(\sigma_t)_{t\ge 0}$ are locally bounded. The volatility process is strictly positive, $\inf_{t\in[0,1]}\sigma_t>0$, $\mathbb{P}^X$-almost surely. 
For all $0\leq t+s\leq1$, $t\ge 0$, $s\ge 0$, with some constants $ C_{\sigma}>0$, and $\alpha>0$, it holds that 
\begin{align}
\label{vola}\E\big[(\sigma_{t+s}-\sigma_{t})^2\big]  \le C_{\sigma} s^{2\alpha}\,.\end{align}
\end{assumption}
Condition \eqref{vola} introduces a regularity parameter $\alpha$, governing the smoothness of the volatility process. The parameter $\alpha$ is crucial, since it will naturally influence convergence rates of spot volatility estimation. Inequality \eqref{vola} is less restrictive than $\alpha$-Hölder continuity, since it does not rule out volatility jumps. For instance, any compound Poisson jump process with a jump size distribution having finite second moments, satisfies \eqref{vola} with $\alpha=1/2$. Since second moments in \eqref{vola} of such a process are bounded by a constant times $(s^2+s)$, i.e.\ the second moment of a Poisson distribution with parameter $s$, this readily follows. Similar bounds for more general jump processes are given, for instance, in Section 13 of \cite{JP}. This is important as empirical evidence for volatility jumps, in particular simultaneous price and volatility jumps, has been reported for intraday high-frequency financial data, see, for instance, \cite{voljumps} and \cite{bibwinkellev}. The presented theory is moreover for general stochastic volatilities, allowing as well for rough volatility. Rough fractional stochastic volatility models recently became popular and are used, for instance, in the macroscopic model of \cite{mathieu2} and \cite{mathieu1}. 

The jump component of \eqref{sm} is illustrated as in \cite{JP} and related literature, where the predictable function $\delta$ is defined on $\Omega\times \mathbb{R}_+\times \mathbb{R}$, and the Poisson random measure $\mu$ is compensated by $\nu(ds,dz)=\lambda(dz)\otimes ds$, with a $\sigma$-finite measure $\lambda$. We impose the following standard condition with a generalized Blumenthal-Getoor or jump activity index $r,\, 0\le r\le 2$.
\begin{assumption}\label{jumps}
Assume that \(\sup_{\omega,x}|\delta(t,x)|/\gamma(x)\) is locally bounded with a non-negative, deterministic function $\gamma$ which satisfies
\begin{align}\label{BG}\int_{\mathbb{R}}(\gamma^r(x)\wedge 1)\lambda(dx)<\infty\,.\end{align}
\end{assumption}
We use the notation $a\wedge b=\min(a,b)$, and $a\vee b=\max(a,b)$, throughout this paper. The assumption is most restrictive in the case $r=0$, when jumps are of finite activity. The larger $r$, the more general jump components are allowed. We will develop results under mild restrictions on $r$.

The process $(X_t)$, which can be decomposed 
\begin{align}\label{smdecomp}X_t=C_t+J_t\,,\end{align}
with the continuous component $(C_t)$, and the c\`{a}dl\`{a}g jump component $(J_t)$, provides a model for the latent efficient log-price process in continuous time.

High-frequency (best) ask order prices from a limit order book at times $t_i^n,~{0\le i\le n}$, on the fix time interval $[0,1]$, cannot be adequately modelled by discrete recordings of $(X_t)$. Instead, we propose the additive model with lower-bounded, one-sided microstructure noise:
\begin{align}\label{lomn}Y_i=X_{t_i^n}+\epsilon_i\,,\,i=0,\ldots,n,~~\epsilon_i\stackrel{iid}{\sim}F_{\eta},\epsilon_i\ge 0\,.\end{align}
The crucial property of the model is that the support of the noise is lower bounded. It is not that important, that this boundary is zero, it could be as well a different constant, or even a regularly varying function over time. The presented methods and results are robust with respect to such model generalizations. We set the bound equal to zero which appears to be the most natural choice for limit orders.
\begin{assumption}\label{noise_ass}
The i.i.d.\ noise $(\epsilon_i)_{0\le i\le n}$, has a cumulative distribution function (cdf) $F_{\eta}$ satisfying
\begin{align}\label{noise_dist}
F_{\eta}(x) =\eta x\big(1+\KLEINO(1)\big) ,\;\mbox{as}~x\downarrow 0\,.
\end{align}
\end{assumption}
This is a nonparametric model in that the extreme value index is $-1$ for the \emph{minimum} domain of attraction close to the boundary. This standard assumption on one-sided noise has been already used by \cite{bound3} and \cite{bound2} within different frameworks. We do not require assumptions about the maximum domain of attraction, moments and the tails of the noise distribution. Parametric examples which satisfy \eqref{noise_dist} are, for instance, the uniform distribution on some interval, the exponential distribution and the standard Pareto distribution with heavy tails. 

The i.i.d.\ assumption on the noise is crucial and generalizations to weakly dependent noise will require considerable work and new proof concepts. Heterogeneity instead, that is, a time-dependent noise level $\eta(t)$, could be included in our asymptotic analysis under mild assumptions.

\section{Construction of spot volatility estimators\label{sec:3}}

We partition the observation interval $[0,1]$ in $h_n^{-1}$ equispaced blocks, $h_n^{-1}\in\N$, and take local minima on each block. We hence obtain for $k=0,\ldots,h_n^{-1}-1$, the local, block-wise minima
\begin{align}\label{localmin}m_{k,n}=\min_{i\in\mathcal{I}_k^n}Y_i~,~\mathcal{I}_k^n=\{i\in\{0,\ldots,n\}: ~t_i^n \in [kh_n,(k+1)h_n)\}\,.\end{align}
While $h_n^{-1}$ is an integer, $nh_n$ is in general not integer-valued. For a simple interpretation, however, one can think of $nh_n$ as an integer-valued sequence which gives the number of noisy observations per block in case of equidistant observations. A spot volatility estimator could be obtained as a localized version of the estimator from Eq.\ (2.9) in \cite{BJR} for the integrated volatility in the analogous model. The idea is that differences $m_{k,n}-m_{k-1,n}$ of local minima estimate differences of efficient prices and a sum of squared differences can be used to estimate the volatility. However, things are not that simple. To determine the expectation of squared differences of local minima we introduce the function
\begin{align}\label{psi}\Psi_n(\sigma^2)=\frac{\pi}{2(\pi-2)} h_n^{-1}\E\Big[\Big(\min_{i\in\{0,\ldots,nh_n-1\}}\hspace*{-.05cm}\big(\sigma B_{\frac{i}{n}}+\epsilon_i\big)-\hspace*{-.05cm}\min_{i\in\{1,\ldots,nh_n\}}\hspace*{-.05cm}\big(\sigma \tilde B_{\frac{i}{n}}+\epsilon_i\big)\Big)^2\Big],\end{align}
where $(B_t)$ and $(\tilde B_t)$ denote two independent standard Brownian motions. In \cite{BJR}, the block length balanced the order of block-wise minimal errors, $(nh_n)^{-1}$ under \eqref{noise_dist}, and the order $h_n^{1/2}$ of the movement of the stochastic semimartingale boundary over a block. For $h_n n^{2/3}\to \infty$, $\Psi_n$ tends to the identity function, such we have that
\begin{align}\label{psiapprox}\Psi_n(\sigma^2)=\sigma^2+\KLEINO(1),\;\mbox{as}~n\to \infty\,.\end{align}
In this asymptotic regime local minima are mainly determined by local minima of the boundary process, such that the first-order approximation equals \eqref{psi} when neglecting the noise $(\epsilon_i)$ on the right-hand side. The half-normal distribution of the minimum of a Brownian motion over an interval and its moments then readily yield \eqref{psiapprox}. A formal proof of \eqref{psiapprox} is contained in Step 3 of the proof of Theorem \ref{propvola} in Section \ref{sec:volacont}. Note that we defined $\Psi_n$ different compared to \cite{BJR}, e.g.\ in their Eq.\ (A.35), with the additional factor $\pi/(\pi-2)$. By the simple asymptotic approximation \eqref{psiapprox}, we do not require $\Psi_n^{-1}$ for a consistent estimator.  

When there are no price jumps, a simple consistent estimator for the spot squared volatility $\sigma_{\tau}^2$ is given by 
\begin{align}\label{simpleestimator}\hat\sigma^2_{\tau-}=\frac{\pi}{2(\pi-2)K_n}\sum_{k=(\lfloor h_n^{-1}\tau\rfloor-K_n)\vee 1}^{\lfloor h_n^{-1}\tau\rfloor-1}h_n^{-1}\big(m_{k,n}-m_{k-1,n})^2\,,\end{align}
for suitable sequences $h_n\to 0$ and $K_n\to\infty$. Using only observations before time $\tau$, the estimator is available on-line at time $\tau\in(0,1]$ during a trading day. For $\tau$ close to 0, when $\lfloor h_n^{-1}\tau\rfloor\le K_n$, the factor $K_n^{-1}$ can be adjusted to get an average. Since this is unimportant for asymptotic theory, we keep $K_n^{-1}$ for a simple notation. Working with ex-post data over the whole interval, instead of using only observations before time $\tau$, one may use as well
\begin{align}\label{simpleestimator2}\hat\sigma^2_{\tau+}=\frac{\pi}{2(\pi-2)K_n}\sum_{k=\lfloor h_n^{-1}\tau\rfloor+1}^{(\lfloor h_n^{-1}\tau\rfloor+K_n)\wedge (h_n^{-1}-1)}h_n^{-1}\big(m_{k,n}-m_{k-1,n})^2\,,\end{align}
or an estimator with an average centred around time $\tau\in(0,1)$. The difference between the two estimators \eqref{simpleestimator2} and \eqref{simpleestimator} can be used to infer a possible jump in the volatility process at time $\tau\in(0,1)$, similarly as in \cite{bibwinkellev}. 

To construct confidence intervals for the spot volatility, it is useful to establish also a spot quarticity estimator:
\begin{align}\label{quartestimator}\widehat{{\sigma^4_{\tau}}}_-=\frac{\pi}{4(3\pi-8)K_n}\sum_{k=(\lfloor h_n^{-1}\tau\rfloor-K_n)\vee 1}^{\lfloor h_n^{-1}\tau\rfloor-1}h_n^{-2}\big(m_{k,n}-m_{k-1,n})^4\,.\end{align}

A spot volatility estimator which is robust with respect to jumps in $(X_t)$ is obtained with threshold versions of these estimators. We truncate differences of local minima whose absolute values exceed a threshold $u_n= \beta\cdot h_n^{\kappa},\kappa\in(0,1/2)$, with some positive constant $\beta$, which leads to
\begin{align}\label{truncatedestimator}\hat\sigma^{2,(tr)}_{\tau-}=\frac{\pi}{2(\pi-2)K_n}\sum_{k=(\lfloor h_n^{-1}\tau\rfloor-K_n)\vee 1}^{\lfloor h_n^{-1}\tau\rfloor-1}h_n^{-1}\big(m_{k,n}-m_{k-1,n})^2\1_{\{|m_{k,n}-m_{k-1,n}|\le u_n\}}\,,\end{align}
and analogous versions of the estimators \eqref{simpleestimator2} and \eqref{quartestimator}.

\section{Asymptotic results\label{sec:4}}
We establish asymptotic results for equidistant observations, $t_i^n=i/n$. We begin with the asymptotic theory in a setup without jumps in $(X_t)$. 
\begin{theorem}[Stable central limit theorem for continuous $(X_t)$]\label{propvola}
Set $h_n$, such that $h_n n^{2/3}\to \infty$, 
and $K_n=C_K h_n^{\delta -2\alpha/(1+2\alpha)}$ for arbitrary $\delta$, $0<\delta<2\alpha/(1+2\alpha)$, and with some constant $C_K>0$. If $(X_t)$ is continuous, i.e.\ $J_t=0$ in \eqref{smdecomp}, under Assumptions \ref{sigma} and \ref{noise_ass}, the spot volatility estimator \eqref{simpleestimator} is consistent, $\hat\sigma^2_{\tau-}\stackrel{\P}{\rightarrow} \sigma_{\tau-}^2$, and satisfies the stable central limit theorem
\begin{align}\label{spotclt}K_n^{1/2}\Big(\hat\sigma^2_{\tau-}-\Psi_n\big(\sigma_{\tau-}^2\big)\Big) \stackrel{st}{\longrightarrow} \mathcal{N}\Big(0,\frac{7\pi^2/4-2\pi/3-12}{(\pi-2)^2}\sigma^4_{\tau-}\Big)\,.\end{align}
\end{theorem}
There is only a difference between $\sigma_{\tau}^2$ and its left limit $\sigma_{\tau-}^2$ in case of a volatility jump at time $\tau$. In particular, the estimator is as well consistent for $\sigma_{\tau}^2$, for any fix $\tau\in(0,1)$. The convergence rate $K_n^{-1/2}$ gets arbitrarily close to $n^{ -2\alpha/(3+6\alpha)}$, which is optimal in our model. The optimal rate is according to \cite{BJR} attained for $h_n\propto n^{-2/3}$, and for $K_n\propto h_n^{-2\alpha/(1+2\alpha)}$, i.e.\ $\delta\downarrow 0$. In the important special case when $ \alpha=1/2$, for a semimartingale volatility, the rate is arbitrarily close to $n^{ -1/6}$. This is much faster than the optimal rate of convergence in the model with additive centred microstructure noise, which is known to be $n^{ -1/8}$, see \cite{hoffmann} and \cite{hautsch}. The constant in the asymptotic variance is obtained from several variance and covariance terms including (squared) local minima and is approximately 2{.}44. The function $\Psi_n$ was shown to be monotone and invertible in \cite{BJR} and $\Psi_n$ and its inverse $\Psi_n^{-1}$ can be approximated using Monte Carlo simulations, see Section \ref{sec:5.1}. The asymptotic distribution of the estimator does not hinge on the noise level $\eta$, different to methods for centred noise. Hence, we do not require any pre-estimation of noise parameters and the theory directly extends to a time-varying noise level $\eta(t)$ in \eqref{noise_dist} under the mild assumption that $0<\eta(t)<\infty$, for all $t$. The stable convergence in \eqref{spotclt} is stronger than weak convergence and is important, since the limit distribution is \emph{mixed} normal depending on the stochastic volatility. We refer to \cite{JP}, Section 2.2.1, for an introduction to stable convergence. For a normalized central limit theorem, we can use the spot quarticity estimator  \eqref{quartestimator}.
\begin{proposition}[Feasible central limit theorem]\label{propvola2}Under the conditions of Theorem \ref{propvola}, the spot quarticity estimator \eqref{quartestimator} 
is consistent, such that we get for the spot volatility estimation the normalized central limit theorem
\begin{align}K_n^{1/2}\frac{\pi-2}{\sqrt{\widehat{{\sigma^4_{\tau}}}_-(7\pi^2/4-2\pi/3-12)}}\Big(\hat\sigma^2_{\tau-}-\Psi_n\big(\sigma_{\tau-}^2\big)\Big) \stackrel{d}{\longrightarrow} \mathcal{N}(0,1)\,.\end{align}
\end{proposition}
Proposition \ref{propvola2} yields asymptotic confidence intervals for spot volatility estimation. For $q\in (0,1)$, it holds true that
\begin{align*}&\P\bigg(\sigma_{\tau-}^2\in\bigg[\Psi_n^{-1}\Big(\hat\sigma^2_{\tau-}-\frac{\sqrt{\widehat{{\sigma^4_{\tau}}}_-(7\pi^2/4-2\pi/3-12)}}{\pi-2}K_n^{-1/2}\Phi^{-1}(1-q/2)\Big)\,,\\
&\quad\quad\quad \Psi_n^{-1}\Big(\hat\sigma^2_{\tau-}+\frac{\sqrt{\widehat{{\sigma^4_{\tau}}}_-(7\pi^2/4-2\pi/3-12)}}{\pi-2}K_n^{-1/2}\Phi^{-1}(1-q/2)\Big)\bigg]\bigg)\to 1-q\,,\end{align*}
by monotonicity of $\Psi_n^{-1}$, with $\Phi$ the cdf of the standard normal distribution. Since $\Psi_n^{-1}$ is differentiable by the result from Eq.\ (A.35) in \cite{BJR} and the derivative is $\big(\Psi_n^{-1}\big)'=1+\KLEINO(1)$ by \eqref{psiapprox}, the delta method (for stable convergence) yields as well asymptotic confidence intervals and the central limit theorem 
\begin{align}\label{spotclt2}K_n^{1/2}\Big(\Psi_n^{-1}\big(\hat\sigma^2_{\tau-}\big)-\sigma_{\tau-}^2\Big) \stackrel{st}{\longrightarrow} \mathcal{N}\Big(0,\frac{7\pi^2/4-2\pi/3-12}{(\pi-2)^2}\sigma^4_{\tau-}\Big)\,.\end{align}
We may not simply replace $\Psi_n\big(\sigma_{\tau-}^2\big)$ in \eqref{spotclt} by its first-order approximation $\sigma_{\tau-}^2$, or $\Psi_n^{-1}\big(\hat\sigma^2_{\tau-}\big)$ in \eqref{spotclt2} by $\hat\sigma^2_{\tau-}$, since the biases do not converge to zero sufficiently fast. That is, $(\hat\sigma^2_{\tau-}-\sigma_{\tau-}^2) =\mathcal{O}_{\P}\big(K_n^{-1/2}\big)$ does in general not hold true. Furthermore, if the condition $h_n n^{2/3}\to\infty$ is violated, the central limit theorems do not apply.
\begin{theorem}[Stable central limit theorem with jumps in $(X_t)$]\label{propvola3}
Set $h_n$, such that $h_n n^{2/3}\to \infty$, 
and $K_n=C_K h_n^{\delta -2\alpha/(1+2\alpha)}$ for arbitrary $\delta$, $0<\delta<2\alpha/(1+2\alpha)$, and with some constant $C_K>0$. Under Assumptions \ref{sigma}, \ref{jumps} and \ref{noise_ass}, with 
\begin{align}r< \frac{2+2\alpha}{1+2\alpha}\,,\end{align}
the truncated spot volatility estimator \eqref{truncatedestimator} with
\begin{align}\kappa\in \Big(\frac{1}{2-r}\frac{\alpha}{2\alpha+1},\frac12\Big)\,,\end{align} 
is consistent, $\hat\sigma^{2,(tr)}_{\tau-}\stackrel{\P}{\rightarrow} \sigma_{\tau-}^2$, and satisfies the stable central limit theorem
\begin{align}\label{spotclttr}K_n^{1/2}\Big(\hat\sigma^{2,(tr)}_{\tau-}-\Psi_n\big(\sigma_{\tau-}^2\big)\Big) \stackrel{st}{\longrightarrow} \mathcal{N}\Big(0,\frac{7\pi^2/4-2\pi/3-12}{(\pi-2)^2}\sigma^4_{\tau-}\Big)\,.\end{align}
\end{theorem}
In order to obtain a central limit theorem at (almost) optimal rate, we thus have to impose mild restrictions on the jump activity. For the standard model with a semimartingale volatility, i.e.\ $\alpha=1/2$, we need that $r<3/2$, and for $\alpha=1$ we have the stronger condition that $r<4/3$. These conditions are equivalent to the ones of Theorem 1 in \cite{bibwinkel2}, which gives a central limit theorem for spot volatility estimation under similar assumptions on $(X_t)$, but with slower rate of convergence for centred microstructure noise. Using a truncated quarticity estimator with the same thresholding yields again a feasible central limit theorem and asymptotic confidence intervals.
\begin{remark}
From a theoretical point of view one might ponder why we do not work out an asymptotic theory for $h_n\propto n^{-2/3}$, when noise and efficient price both influence the asymptotic distribution of the local minima. However, in this balanced case, the asymptotic distribution is infeasible. For this reason, \cite{BJR} could not establish a central limit theorem for their integrated volatility estimator. Moreover, their estimator was only implicitly defined depending on the unknown function $\Psi_n^{-1}$. Even imposing a parametric assumption on the noise as an exponential distribution would not render a feasible limit theory for $h_n\propto n^{-2/3}$, see the discussion in \cite{BJR}. 
Choosing $h_n$, such that $h_n n^{2/3}\to\infty$ slowly, yields instead a simple, explicit and consistent estimator and a feasible central limit theorem for spot volatility estimation. In particular, we use $\Psi_n$ only for the bias-correction of the simple estimator, while the estimator itself and the (estimated) asymptotic variance do not hinge on $\Psi_n$. Central limit theorems for spot volatility estimators are in general only available at almost optimal rates, when the variance dominates the squared bias in the mean squared error, see, for instance, Theorem 13.3.3 and the remarks below in \cite{JP}. Therefore, \eqref{spotclt} is the best achievable central limit theorem. Our choice of $h_n$ avoids moreover strong assumptions on the noise that would be inevitable for smaller blocks. Our numerical work will demonstrate that the presented asymptotic results are useful in practice and facilitate efficient inference on the spot volatility. In particular, Section \ref{sec:5.2} evolves around the question how to choose block lengths in practice. 
\end{remark}

\section{Implementation and simulations\label{sec:5}} 
\subsection[Monte Carlo approximation of Psi]{Monte Carlo approximation of $\Psi_n$\label{sec:5.1}}
Although the function $\Psi_n$ from \eqref{psi} tends to the identity asymptotically, it has a crucial role for a bias correction of our estimator in \eqref{spotclt}. We can compute the function numerically based on a Monte Carlo simulation. Hence, we have to compute $\Psi_n(\sigma^2)$ as a Monte Carlo mean over many iterations and over a fine grid of values for the squared volatility. Then, we can also numerically invert the function and use $\Psi_n^{-1}(\,\cdot\,)$. To make this procedure feasible without too high computational expense we require an algorithm to efficiently sample from the law of the local minima for some given $n$ and block length $h_n$.

Consider for $nh_n\in\N$, with $Z_i\stackrel{iid}{\sim}\mathcal{N}(0,1)$, and the observation errors $(\epsilon_k)_{k\ge 0}$, the minimum
\[M_1^{nh_n}:=\min_{k=1,\ldots,nh_n}\Big(\frac{\sigma}{\sqrt{n}}\sum_{i=1}^k Z_i+\epsilon_k\Big)\,,\]
for some fix $\sigma>0$, and for $l\in\{0,\ldots,nh_n\}$:
\[M_l^{nh_n}:=\min_{k=l,\ldots,nh_n}\Big(\frac{\sigma}{\sqrt{n}}\sum_{i=0}^k Z_i+\epsilon_k\Big)\,,\]
where we set $Z_0:=0$. Since
\[\Psi_n(\sigma^2)=\frac12\frac{\pi}{\pi-2}h_n^{-1}\E\Big[\big(M_0^{nh_n-1}-M_1^{nh_n}\big)^2\Big]\,,\]
with $M_0^{nh_n-1}$ generated independently from $M_1^{nh_n}$, we want to simulate samples distributed as $M_0^{nh_n-1}$ and $M_1^{nh_n}$, respectively. Note that for finite $nh_n$, there is not an exact equality between the moments of $M_0^{nh_n-1}$ and $M_1^{nh_n}$, what can be relevant  in particular for moderate values of $nh_n$. As in the simulation of Section \ref{sec:5.2}, we implement exponentially distributed observation errors $(\epsilon_k)$, with some given noise level $\eta$. In data applications, we can do the same with an estimated noise level
\begin{align*}\hat\eta=\bigg(\frac{1}{2n}\sum_{i=1}^n\big(Y_i-Y_{i-1}\big)^2\bigg)^{-1/2}=\eta+\mathcal{O}_{\P}\big(n^{-1/2}\big)\,.\end{align*}
This estimator works for all noise distributions with finite fourth moments. In view of the discussion of the model in \cite{BJR}, exponentially distributed noise is the most natural example satisfying \eqref{noise_dist}. Simulations with other noise distributions led to similar results. This is expected, since the estimator only hinges on local minima and their distribution is asymptotically more determined by the Brownian motion than by the noise distribution. To simulate the local minima for given $n$, $h_n$, $\eta$, and squared volatility $\sigma^2$, in an efficient way we use a specific dynamic programming principle. Observe that
\begin{align*}
M_1^{nh_n}&=\frac{\sigma}{\sqrt{n}}Z_1+\min\Big(\epsilon_1,M_2^{nh_n}\Big)\\
&=\frac{\sigma}{\sqrt{n}}Z_1+\min\Big(\epsilon_1,\frac{\sigma}{\sqrt{n}}Z_2+\min\Big(\epsilon_2,M_3^{nh_n}\Big)\Big)\\
& =\frac{\sigma}{\sqrt{n}}Z_1+\min\bigg( \ldots~\min\bigg(\epsilon_{nh_n-2},\frac{\sigma}{\sqrt{n}}Z_{nh_n-1}+\min\Big(\epsilon_{nh_n-1},\frac{\sigma}{\sqrt{n}}Z_{nh_n}+\epsilon_{nh_n}\Big)\bigg)\ldots\bigg).
\end{align*}
In the baseline noise model, $\epsilon_k\stackrel{iid}{\sim}\text{Exp}(\eta)$, the random variable $\frac{\sigma}{\sqrt{n}}Z_{nh_n}+\epsilon_{nh_n}$ has an exponentially modified Gaussian (EMG) distribution. With any fixed noise distribution, we can easily generate realizations from this convolution. A pseudo random variable which is distributed as $M_1^{nh_n}$ is now generated following the last transformation in the reverse direction. In pseudo code, this reads
\begin{verbatim}
1. Generate U_{nh_n}~ EMG(sigma^2/n,eta)~ Exp(eta)+sigma/sqrt(n)*Norm(1)
2. U_{nh_n-1}=min(U_{nh_n},Exp(eta))+sigma/sqrt(n)*Norm(1)
3. iterate until U_1
\end{verbatim}
where the end point $U_1$ has the target distribution of $M_1^{nh_n}$. In each iteration step, we thus take the minimum of the current state of the process with one independent exponentially distributed random variable and the convolution with one independent normally distributed random variable. To sample from the distribution of $M_0^{nh_n-1}$ instead, we use the same algorithm and just drop the convolution with the normal distribution in the last step. 

This algorithm facilitates a many times faster sampling from the distribution of local minima and numerical approximation 
of $\Psi_n$, compared to running for each value a standard Monte Carlo simulation in that local minima are computed over blocks of length $h_n$. 

\begin{figure}[t]
\begin{center}
\includegraphics[width=6.5cm]{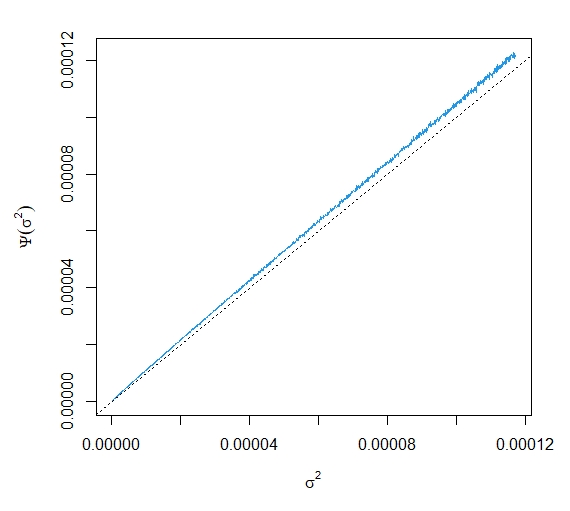}\includegraphics[width=6.5cm]{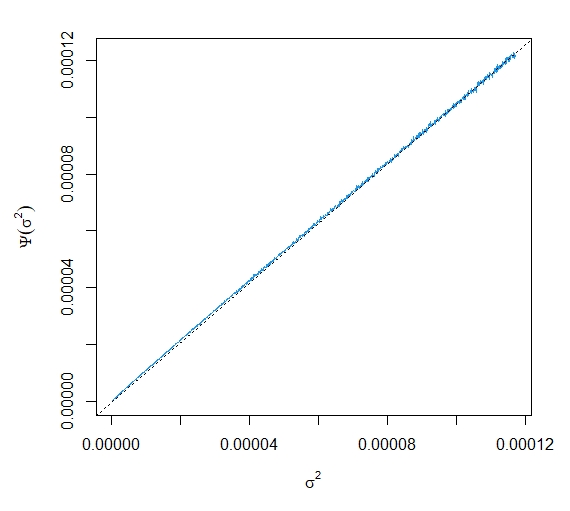}
\caption{Monte Carlo means to estimate $\Psi_n(\sigma^2)$ over a fine grid (blue line) for $n=23{,}400$ and $n\cdot h_n=15$. Left, the dotted line shows the identity function, right the dotted line is a linear function with slope 1{.}046.}\label{Fig:psi}
\end{center}
\end{figure}

Figure \ref{Fig:psi} plots the result of the Monte Carlo approximation of $\Psi_n(\sigma^2)$ for $n=23{,}400$ and $n\cdot h_n=15$, on a grid of 1500 values of $\sigma^2$. In this case, $h_n$ is quite small, but this configuration turns out to be useful below in Section \ref{sec:5.2}. We know that $\Psi_n(\sigma^2)$ is monotone, such that the oscillation of the blue line in Figure \ref{Fig:psi} is due to the inaccuracy of the Monte Carlo means although we use $N=100{,}000$ iterations for each grid point. Nevertheless, we can see that the function is rather close to a linear function with slope $1{.}046$ based on a least squares estimate. The left plot of Figure \ref{Fig:psi} draws a comparison to the identity function which is illustrated by the dotted line, while the plot right-hand side draws a comparison to the linear function with slope $1{.}046$. We see that it is crucial to correct for the bias in \eqref{spotclt} when using such small values of $h_n$. Although the function $\Psi_n(\sigma^2)$ is not exactly linear, a simple bias correction dividing estimates by $1{.}046$ is almost as good as using the more precise numerical inversion based on the Monte Carlo approximation. Since the Monte Carlo approximations of $\Psi_n(\sigma^2)$ look close to linear functions in all considered cases, we report the estimated slopes based on least squares and $N=100{,}000$ Monte Carlo iterations for different values of $h_n$ in Table \ref{tab:1} to summarize concisely how far the distance between the function $\Psi_n(\sigma^2)$ and the identity is. Simulating all iterations for all grid points with our algorithm takes only a few hours with a standard computer. 

\begin{table}[t]
\caption{Regression slopes to measure the bias of estimator \eqref{simpleestimator} and deviation $\Psi_n(\sigma^2)-\sigma^2$.\label{tab:1}}
\begin{center}
\begin{tabular}{|c|c|c|c|c|c|c|}
\hline
 $n\cdot h_n$& 10&15 &25&39& 78& 234\\
\hline
$h_n^{-1}$& 2340& 1560 &936&600&300 &100\\ 
\hline
$h_n\cdot n^{2/3}$& 0{.}350& 0{.}524& 0{.}874&1{.}36 &2{.}73&8{.}18\\
\hline
slope& 1{.}077 &1{.}046&1{.}025&1{.}016 & 1{.}008 &1{.}003\\
\hline
approx.\ bias& 7{.7}\%& 4{.6}\%&2{.5}\%&1{.}6\%&0{.}8\%&0{.}3\%\\
\hline
\end{tabular}
\end{center}
\end{table}

\subsection{Simulation study of estimators \label{sec:5.2}}
\renewcommand\thefootnote{\arabic{footnote}}
\setcounter{footnote}{0}

We simulate $n=23{,}400$ observations corresponding to one observation per second over a (NASDAQ) trading day of 6.5 hours. The efficient price process is simulated from the model
\begin{align*}
dX_t&=\nu_t\sigma_t\,dW_t\,,\\
d\sigma_t^2&= 0{.}0162\cdot \big(0{.}8465 - \sigma_t^2\big) \,dt + 0{.}117 \cdot \sigma_t \, dB_t\,,\\
\nu_t &= \big(6 - \sin(3\pi t/4)\big) \cdot 0{.}002\,,~t\in[0,1]\,.
\end{align*}
The factor $(\nu_t)$ generates a typical U-shaped intraday volatility pattern. $(W_t,B_t)$ is a two-dimensional Brownian motion with leverage $d[W,B]_t=-0{.}2\,dt$. The stochastic volatility component has several realistic features and the simulated model is in line with recent literature, see \cite{bibwinkellev} and references therein. We do not include a drift in $X_t$ to avoid introducing another process or more parameters. Any drift evolving within a reasonable range of values will not affect the presented numerical results. Observations with lower-bounded, one-sided microstructure noise are generated by
\begin{align*}Y_i=X_{\frac{i}{n}}+\epsilon_i\,,~0\le i\le n\,,
\end{align*}
with exponentially distributed noise, $\epsilon_i\stackrel{iid}{\sim}\text{Exp}(\eta)$, with $\eta=10{,}000$. The noise variance is then rather small, but this is in line with stylized facts of real NASDAQ data as, for instance, those analysed in \cite{hautsch} and \cite{bibwinkellev}.\footnote{Note that the noise level estimate is analogous to the one used for regular market microstructure noise. Typical noise levels obtained for trades of e.g.\ Apple are approx.\ 15{,}000 and approx.\ 4{,}000 for 3M, see the supplement of \cite{hautsch}. For mid quotes or best ask/bid prices the levels are only slightly larger (variance smaller).}

\begin{figure}[t]
\begin{center}
\includegraphics[width=13cm]{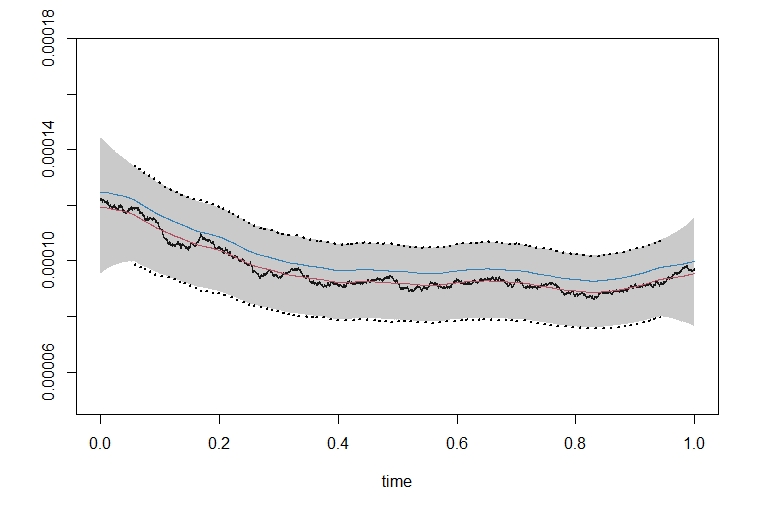}
\caption{True and estimated spot volatility with pointwise confidence sets.}\label{Fig:volaest}
\end{center}
\end{figure}

The black line in Figure \ref{Fig:volaest} shows a fixed path of the squared volatility. We fix this path for the following Monte Carlo simulation and generate new observations of $(X_t)$ and $(Y_i)$ in each iteration according to our model. The blue line in Figure \ref{Fig:volaest} gives the estimated volatility by the Monte Carlo means over $N=50{,}000$ iterations based on $n\cdot h_n=15$ observations per block using the non-adjusted estimator \eqref{simpleestimator}, but with windows which are centred around the block on that we estimate the spot volatility, i.e.\ windows centred around the time $\tau$, and with $K_n=180$. We plot estimates on each block, where the estimates close to the boundaries rely on less observations. The red line gives the bias-corrected volatility estimates using the numerically evaluated function $\Psi_n$, based on the algorithm from Section \ref{sec:5.1} with $n\cdot h_n=15$ and $n=23{,}400$. We determined the values $n\cdot h_n=15$ and $K_n=180$ as suitable values to obtain a small mean squared error. In fact, the choice of $K_n=180$ is rather large in favour of a smaller variance what yields a rather smooth estimated spot volatility in Figure \ref{Fig:volaest}. The estimated volatility hence appears smoother compared to the true semimartingale volatility, but the intraday pattern is well captured by our estimation. We expect that this is typically an appealing implementation in practice as smaller $K_n$ results in a larger variance. Choosing $K_n=180$ rather large, we have to use quite small block sizes $h_n$, to control the overall bias of the estimation. Since $h_n\cdot n^{2/3}\approx 0{.}52$ is small, the bias correction becomes crucial here. Still, our asymptotic results work well for this implementation. This can be seen by the comparison of pointwise empirical 10\% and 90\% quantiles from the Monte Carlo iterations illustrated by the grey area and the 10\% and 90\% quantiles of the limit normal distribution with the asymptotic variance from \eqref{spotclt}. The latter are drawn as dotted lines for the blocks with larger distance than $K_n/2$ from the boundaries where the variances are of order $K_n^{-1}$. Close to the boundaries the empirical variances increase due to the smaller number of blocks used for the estimates. Moreover, the bias correction which is almost identical to dividing each estimate by $1{.}046$, correctly scales the simple estimates which have a significant positive bias for the chosen tuning parameters. Overall, our asymptotic results provide a good finite-sample fit even though we have $h_n\cdot n^{2/3}<1$ here. Note, however, that $\sigma_t \cdot \eta\approx 100$, and our asymptotic expansion requires in fact that $h_n^{3/2} n\sigma_t \eta$ is large when taking constants into account. Since the simulated scenario uses realistic values, we recommend similar block lengths for applications to real high-frequency financial data. According to the summary statistics in the supplement of \cite{hautsch}, some assets exhibit higher noise-to-signal ratios and for those larger blocks are preferable.

\begin{table}[t]
\caption{Summary statistics of estimation for different values of $h_n$ and $K_n$, MSD = mean standard deviation, MAB = mean absolute bias, MABC = MAB of bias-corrected estimator.\label{tab:2}}
\begin{center}
\begin{tabular}{|c|c|c|c|c|c|c|c|c|c|}
\hline
 \rowcolor{lightgray}$K_n$&\multicolumn{3}{|c|}{120}&\multicolumn{3}{|c|}{180}&\multicolumn{3}{|c|}{240}\\
\hline
\cellcolor{gray!20}$nh_n$&MSD&MAB&MABC&MSD&MAB&MABC&MSD&MAB&MABC\\
\hline
\cellcolor{gray!20}10& 14.6&7.59&0.73&12.0&7.51&0.90&10.5&7.60&1.13\\ 
\hline
\cellcolor{gray!20}15& 14.4&4.59&0.88&11.8&4.57&1.17&10.3&4.46&1.43\\
\hline
\cellcolor{gray!20}25& 14.3&2.56&1.24&11.8&2.63&1.66&10.3&2.86&1.91\\
\hline
\cellcolor{gray!20}78&14.7&2.44&2.52&12.3&3.53&3.42&11.0&4.33&4.16 \\
\hline
\end{tabular}
\end{center}
\footnotesize{All values multiplied with factor $10^6$.}
\end{table}

Table \ref{tab:2} summarizes the performance of the estimation along different choices of $nh_n$ and $K_n$. We give the following quantities:
\begin{enumerate}
\item MSD: the \textbf{m}ean \textbf{s}tandard \textbf{d}eviation of $N$ iterations averaged over all grid points; 
\item MAB: the \textbf{m}ean \textbf{a}bsolute \textbf{b}ias of $N$ iterations averaged over all grid points and for estimator \eqref{simpleestimator} without any bias correction; 
\item MABC: the \textbf{m}ean \textbf{a}bsolute \textbf{b}ias of $N$ iterations averaged over all grid points and for estimator \eqref{simpleestimator} with a simple bias \textbf{c}orrection dividing estimates by the factors given in Table \ref{tab:1}.
\end{enumerate} 
All results are based on $N=50{,}000$ Monte Carlo iterations. First of all, the values used for Figure \ref{Fig:volaest} are not unique minimizers of the mean squared error. Several other combinations given in Table \ref{tab:2} render equally well results. Overall, the performance is comparable within a broad range of block lengths and window sizes. The variances decrease for larger $K_n$, while the bias increases with larger $K_n$ for fixed $h_n$. Important for the bias is the total window size, $K_n\cdot h_n$, over that the volatility is approximated constant for the estimation. The variance only depends on $K_n$, changing the block length for fix $K_n$ does not significantly affect the variance. While the MSD is hence almost constant within the columns of Table \ref{tab:2}, the bias after correction, MABC, increases from top down due to the increasing window size. Without the bias correction two effects interfere for MAB. Larger blocks reduce the systematic bias due to $\Psi_n(\sigma_t^2)-\sigma_t^2$, but the increasing bias due to the increasing window size prevails for $n\cdot h_n=78$, and the two larger values of $K_n$.

\section{Proofs\label{sec:6}}
\subsection{Law of the integrated negative part of a Brownian motion}
A crucial lemma for our theory is on an upper bound for the cdf of the integrated negative part of a Brownian motion. We prove a lemma based on a generalization of Lévy's arc-sine law by \cite{takacs1996}. The result is in line with the conjecture in Eq.\ (261) of \cite{Janson} where one finds an expansion of the density with a precise constant of the leading term. Denote by $f_+$ the positive part and by $f_-$ the negative part of some real-valued function $f$.
\begin{lem}\label{lem3}For a standard Brownian motion $(W_t)_{t\ge 0}$, it holds that
\[\P\Big(\int_0^1(W_t)_-\,dt\le x\Big)=\mathcal{O}(x^{1/3}), ~x\to 0\,.\]
\end{lem}
\begin{proof}
Observe the equality in distribution $\int_0^1(W_t)_-\,dt\stackrel{d}{=}\int_0^1(W_t)_+\,dt$, such that
\[\P\Big(\int_0^1(W_t)_-\,dt\le x\Big)=\P\Big(\int_0^1(W_t)_+\,dt\le x\Big)\,,x>0\,.\]
For any $\epsilon>0$, the inequality
\[\int_0^1(W_t)_+\,dt\ge \int_0^1W_t \cdot\1(W_t>\epsilon)\,dt\ge \epsilon \int_0^1\1(W_t>\epsilon)\,dt\]
leads us to
\begin{align*}\P\Big(\int_0^1(W_t)_+\,dt\le x\Big)&\le \P\Big(\epsilon\int_0^1\1(W_t>\epsilon)\,dt\le x\Big)\\
&=\P\Big(1-\int_0^1\1(W_t\le \epsilon)\,dt\le x/\epsilon\Big)\\
&=\P\Big(\int_0^1\1(W_t\le \epsilon)\,dt\ge 1-x/\epsilon\Big)\,.
\end{align*}
Using (15) and (16) from \cite{takacs1996}, we obtain that
\begin{align*}
\P\Big(\int_0^1\1(W_t\le \epsilon)\,dt\ge 1-x/\epsilon\Big)&=\frac{1}{\pi}\int_{1-x/\epsilon}^{1}\frac{\exp(-\epsilon^2/(2u))}{\sqrt{u(1-u)}} du +2\Phi(\epsilon)-1\,,
\end{align*}
with $\Phi$ the cdf of the standard normal distribution. Thereby, we obtain that
\begin{align*}\P\Big(\int_0^1(W_t)_+\,dt\le x\Big)&\le \frac{1}{\pi}\int_{1-x/\epsilon}^{1}\frac{\exp(-\epsilon^2/(2u))}{\sqrt{u(1-u)}} du +2\int_0^{\epsilon}\frac{\exp(-u^2/2)}{\sqrt{2\pi}}du\,,
\end{align*}
and elementary bounds give the upper bound
\begin{align*}\P\Big(\int_0^1(W_t)_+\,dt\le x\Big)&\le \frac{2}{\pi}\sqrt{\frac{x}{\epsilon}}\frac{1}{\sqrt{1-x/\epsilon}}  +\frac{2\epsilon}{\sqrt{2\pi}}\,.
\end{align*}
Choosing $\epsilon=x^{1/3}$, we obtain the upper bound
\begin{align*}\P\big(\int_0^1(W_t)_+\,dt\le x\big)&\le \frac{2}{\pi}x^{1/3}\frac{1}{\sqrt{1-x^{2/3}}}  +\frac{2x^{1/3}}{\sqrt{2\pi}}\,.
\end{align*}
\end{proof}

\subsection{Asymptotics of the spot volatility estimation in the continuous case\label{sec:volacont}}
\subsubsection{Proof of Theorem \ref{propvola}}
In the sequel, we write $A_n\lesssim B_n$ for two real sequences, if there exists some $n_0\in\N$ and a constant $K$, such that $A_n\le K B_n$, for all $n\ge n_0$.\\[.2cm]
{\bf{Step 1}}\\[.1cm]
In the first step, we prove the approximation
\begin{align*}\hat\sigma^2_{\tau-}&=\frac{\pi}{2(\pi-2)K_n}\sum_{k=(\lfloor h_n^{-1}\tau\rfloor-K_n)\vee 1}^{\lfloor h_n^{-1}\tau\rfloor-1}h_n^{-1}\big(m_{k,n}-m_{k-1,n})^2\\
&=\frac{\pi}{2(\pi-2)K_n}\sum_{k=(\lfloor h_n^{-1}\tau\rfloor-K_n)\vee 1}^{\lfloor h_n^{-1}\tau\rfloor-1}h_n^{-1}\big(\tilde m_{k,n}-\tilde m^*_{k-1,n})^2+\mathcal{O}_{\P}\big(h_n^{\alpha\wedge 1/2}\big)
\end{align*}
with 
\begin{align*}
\tilde m_{k,n}&=\min_{i\in\mathcal{I}_k^{n}}\big(\epsilon_i+\sigma_{(k-1)h_n}(W_{t_i^n}-W_{kh_n})\big)\,,~\text{and}\\
\tilde m^*_{k-1,n}&=\min_{i\in\mathcal{I}_{k-1}^{n}}\big(\epsilon_i-\sigma_{(k-1)h_n}(W_{kh_n}-W_{t_i^n})\big)\,.
\end{align*}
We show that for $k\in\{1,\ldots,h_n^{-1}-1\}$, it holds that
\begin{align}\label{eqlem1}m_{k,n}-m_{k-1,n}=\tilde m_{k,n}-\tilde m^*_{k-1,n}+\KLEINO_{\P}\big(h_n^{1/2}\big)\,.\end{align}
We subtract $X_{kh_n}$ from $m_{k,n}$ and $m_{k-1,n}$, and use that it holds for all $i$ that
\[\big(Y_{i}-X_{kh_n}\big)-\big(X_{t_i^n}-\big(X_{kh_n}+\sigma_{(k-1)h_n}(W_{t_i^n}-W_{kh_n})\big)\big)=\big(\sigma_{(k-1)h_n}(W_{t_i^n}-W_{kh_n})+\epsilon_{i}\big)\,.\]
This implies that
\begin{align*}&\min_{i\in\mathcal{I}_{k}^n}\big(Y_{i}-X_{kh_n}\big)-\max_{i\in\mathcal{I}_{k}^n}\big(X_{t_i^n}-\big(X_{kh_n}+\sigma_{(k-1)h_n}(W_{t_i^n}-W_{kh_n})\big)\big)\\
&\quad\le \min_{i\in\mathcal{I}_{k}^n}\big(\sigma_{(k-1)h_n}(W_{t_i^n}-W_{kh_n})+\epsilon_{i}\big)\,.\end{align*}
Changing the roles of $\big(Y_{i}-X_{kh_n}\big)$ and $\big(\sigma_{(k-1)h_n}(W_{t_i^n}-W_{kh_n})+\epsilon_{i}\big)$, we obtain by the analogous inequalities and the triangle inequality, with $M_t:=X_{kh_n}+\int_{kh_n}^{t}\sigma_{(k-1)h_n}\,dW_s$, that
\begin{align*}\Big|m_{k,n}-X_{kh_n}-\tilde m_{k,n}\Big|&\le \max_{i\in\mathcal{I}_{k}^n}\big|X_{t_i^n}-M_{t_i^n}\big|\le \sup_{t\in[kh_n,(k+1)h_n]}\big|X_t-M_t\big|\\
&\le \sup_{t\in[kh_n,(k+1)h_n]}\Big|C_t-C_{kh_n}-\smallint_{kh_n}^t\sigma_{(k-1)h_n}\,dW_s\Big|\,.
\end{align*}
We write $(C_t)$ for $(X_t)$ to emphasize continuity, see \eqref{smdecomp}. \eqref{eqlem1} follows from
\begin{align}\label{help1}\sup_{t\in[kh_n,(k+1)h_n]}\Big|C_t-C_{kh_n}-\smallint_{kh_n}^t\sigma_{(k-1)h_n}\,dW_s\Big|=\KLEINO_{\P}(h_n^{1/2})\,,\end{align}
and the analogous estimate for $m_{k-1,n}$ and $\tilde m^*_{k-1,n}$. We decompose 
\begin{align*}\sup_{t\in[kh_n,(k+1)h_n]}\Big|C_t-C_{kh_n}-\smallint_{kh_n}^t\sigma_{(k-1)h_n}\,dW_s\Big|&\le \sup_{t\in[kh_n,(k+1)h_n]}\Big|\smallint_{kh_n}^t(\sigma_s-\sigma_{(k-1)h_n})\,dW_s\Big|\\
&\quad +\sup_{t\in[kh_n,(k+1)h_n]}\int_{kh_n}^t |a_s|ds\,.\end{align*}
Under Assumption \ref{sigma}, we can assume that $(\sigma_t)$ and $(a_t)$ are bounded on $[0,1]$ by the localization from Section 4.4.1 in \cite{JP}. Using It\^{o}'s isometry and Fubini's theorem, we obtain that
\begin{align*}
\E\Big[\Big(\int_{kh_n}^t(\sigma_s-\sigma_{(k-1)h_n})\,dW_s\Big)^2\Big]&=\E\Big[\int_{kh_n}^t(\sigma_s-\sigma_{(k-1)h_n})^2\,ds\Big]\\
&=\int_{kh_n}^t\E\big[(\sigma_s-\sigma_{(k-1)h_n})^2\big]\,ds\,,\end{align*}
such that Assumption \ref{sigma} yields, for any $t\in[kh_n,(k+1)h_n]$, that
\begin{align*}
\E\Big[\Big(\int_{kh_n}^t(\sigma_s-\sigma_{(k-1)h_n})\,dW_s\Big)^2\Big]&\le C_{\sigma}\int_{kh_n}^t(s-(k-1)h_n)^{2\alpha}\,ds
\\&\le C_{\sigma}(2\alpha+1)^{-1}(t-(k-1)h_n)^{2\alpha+1}=\mathcal{O}\big(h_n^{2\alpha+1}\big)\,.
\end{align*}
By Doob's martingale maximal inequality and since $\sup_{t\in[kh_n,(k+1)h_n]}\int_{kh_n}^t |a_s|ds=\mathcal{O}_{\P}(h_n)$, it holds that
\[\sup_{t\in[kh_n,(k+1)h_n]}\Big|C_t-C_{kh_n}-\smallint_{kh_n}^t\sigma_{(k-1)h_n}\,dW_s\Big|=\mathcal{O}_{\P}\big(h_n^{(1/2+\alpha)\wedge 1}\big)\,.\]
We conclude that \eqref{help1} holds, since $\alpha>0$. Since
\[h_n^{-1}\big(m_{k,n}-m_{k-1,n}\big)\big(m_{k,n}-\tilde m_{k,n}\big)=\mathcal{O}_{\P}\big(h_n^{\alpha\wedge 1/2}\big)\,,\]
and analogously for $(m_{k-1,n}-\tilde m^*_{k-1,n})$, we conclude Step 1 writing
\begin{align*}\big(m_{k,n}-m_{k-1,n}\big)^2-\big(\tilde m_{k,n}-\tilde m^*_{k-1,n}\big)^2&=\big(m_{k,n}-m_{k-1,n}+\tilde m_{k,n}-\tilde m^*_{k-1,n}\big)\\
&\quad\times \big(m_{k,n}-\tilde m_{k,n}+\tilde m^*_{k-1,n}-m_{k-1,n}\big).\end{align*}
{\bf{Step 2}}\\[.1cm]
We bound the bias of the spot volatility estimation using Step 1. For $\lfloor h_n^{-1}\tau\rfloor>K_n$, we obtain from the definition of the function $\Psi_n$ in \eqref{psi} that
\begin{align*}&\E\big[\hat\sigma^2_{\tau-}-\Psi_n\big(\sigma_{\tau-}^2\big)\big]=\frac{\pi}{2(\pi-2)K_n}\sum_{k=(\lfloor h_n^{-1}\tau\rfloor-K_n)\vee 1}^{\lfloor h_n^{-1}\tau\rfloor-1}h_n^{-1}\E\Big[\big( m_{k,n}-m_{k-1,n})^2\Big]-\E\big[\Psi_n\big(\sigma_{\tau-}^2\big)\big]\\
&=\frac{1}{K_n}\frac{\pi}{2(\pi-2)}\sum_{k=(\lfloor h_n^{-1}\tau\rfloor-K_n)\vee 1}^{\lfloor h_n^{-1}\tau\rfloor-1}h_n^{-1}\E\Big[\big(\tilde m_{k,n}-\tilde m^*_{k-1,n})^2\Big]-\E\big[\Psi_n\big(\sigma_{\tau-}^2\big)\big]+\mathcal{O}\big(h_n^{\alpha\wedge 1/2}\big)\\
&=\frac{1}{K_n}\frac{\pi}{2(\pi-2)}\sum_{k=(\lfloor h_n^{-1}\tau\rfloor-K_n)\vee 1}^{\lfloor h_n^{-1}\tau\rfloor-1}\frac{2(\pi-2)}{\pi}\E\big[\Psi_n\big(\sigma_{(k-1)h_n}^2\big)\big]-\E\big[\Psi_n\big(\sigma_{\tau-}^2\big)\big]+\mathcal{O}\big(h_n^{\alpha\wedge 1/2}\big)\\
&\lesssim \frac{1}{K_n}\sum_{k=(\lfloor h_n^{-1}\tau\rfloor-K_n)\vee 1}^{\lfloor h_n^{-1}\tau\rfloor-1}\E\big[\sigma_{(k-1)h_n}^2-\sigma_{\tau-}^2\big]+\mathcal{O}\big(h_n^{\alpha\wedge 1/2}\big)\\
&\lesssim \frac{1}{K_n}\sum_{k=(\lfloor h_n^{-1}\tau\rfloor-K_n)\vee 1}^{\lfloor h_n^{-1}\tau\rfloor-1}\E\big[\sigma_{(k-1)h_n}-\sigma_{\tau-}\big]+\mathcal{O}\big(h_n^{\alpha\wedge 1/2}\big)\\
&\lesssim \frac{1}{K_n}\sum_{k=(\lfloor h_n^{-1}\tau\rfloor-K_n)\vee 1}^{\lfloor h_n^{-1}\tau\rfloor-1}\Big(\E\big[\big(\sigma_{(k-1)h_n}-\sigma_{\tau-}\big)^2\big]\Big)^{1/2}+\mathcal{O}\big(h_n^{\alpha\wedge 1/2}\big)\\
&=\mathcal{O}\big((K_n\,h_n)^{\alpha}\big)=\KLEINO\big(h_n^{\alpha/(1+2\alpha)}\big)=\KLEINO\big(K_n^{-1/2}\big)\,.
\end{align*}
The first $\lesssim$ estimate is in fact an equality up to an additional factor $(1+\KLEINO(1))$, since $\Psi_n'(x)=1+\KLEINO(1)$, for all $x\ge 0$, exploiting the above mentioned differentiability based on Eq.\ (A.35) from \cite{BJR}. For the asymptotic upper bounds we used the binomial formula 
\[\sigma_{(k-1)h_n}^2-\sigma_{\tau-}^2=\big(\sigma_{(k-1)h_n}-\sigma_{\tau-}\big)\big(\sigma_{(k-1)h_n}+\sigma_{\tau-}\big)\le 2\,C\, \big(\sigma_{(k-1)h_n}-\sigma_{\tau-}\big)\,,\]
exploiting as in Step 1 that $(\sigma_t)$ is bounded with some upper bound $C$, and Hölder's inequality to conclude with \eqref{vola} from Assumption \ref{sigma}.
Finally, we used that $(\alpha\wedge 1/2)>\alpha/(2\alpha+1)$ for all $\alpha$.  \\[.2cm]
{\bf{Step 3}}\\[.1cm]
For the consistency of $\hat\sigma^2_{\tau-}$, we prove that
\begin{align}\label{expsimple}\E\Big[\hat\sigma^2_{\tau-}-\sigma_{\tau-}^2\Big]=\KLEINO(1)\,.\end{align}
This includes a proof of \eqref{psiapprox}. Denote by $\P_{\sigma_{(k-1)h_n}}$ the regular conditional probabilities conditioned on $\sigma_{(k-1)h_n}$, and  $\E_{\sigma_{(k-1)h_n}}$ the expectations with respect to the conditional measures. We obtain by the tower rule that
\begin{align}\label{tower}\E\big[h_n^{-1}\big(\tilde m_{k,n}-\tilde m^*_{k-1,n})^2\big]&=\E\big[h_n^{-1}\E_{\sigma_{(k-1)h_n}}\big[\big(\tilde m_{k,n}-\tilde m^*_{k-1,n})^2\big]\big]\\
&\notag=\E\Big[\E_{\sigma_{(k-1)h_n}}\big[\big(h_n^{-1/2}\tilde m_{k,n})^2\big]+\E_{\sigma_{(k-1)h_n}}\big[\big(h_n^{-1/2}\tilde m^*_{k-1,n})^2\big]\\
&\notag\quad -2\,\E_{\sigma_{(k-1)h_n}}\big[h_n^{-1/2}\tilde m_{k,n}\big]\E_{\sigma_{(k-1)h_n}}\big[h_n^{-1/2}\tilde m^*_{k-1,n}\big]\Big]\,,
\end{align}
by the conditional independence of $\tilde m_{k,n}$ and $\tilde m^*_{k-1,n}$. 

We establish and use an approximation of the tail probabilities of $(\tilde m_{k,n})$ and $(\tilde m^*_{k-1,n})$, respectively. For $x\in\R$, we have that
\begin{align}
&\P_{\sigma_{(k-1)h_n}}\Big(h_n^{-1/2}\min_{i\in\mathcal{I}_{k}^n} \big(\epsilon_{i}+\sigma_{(k-1)h_n}(W_{t_i^n}-W_{kh_n})\big)>x\sigma_{(k-1)h_n}\Big)\\
~&\notag=\P_{\sigma_{(k-1)h_n}}\Big(\min_{i\in\mathcal{I}_{k}^n} \big(h_n^{-1/2}\big(W_{t_i^n}-W_{kh_n}\big)+h_n^{-1/2}\sigma_{(k-1)h_n}^{-1}\epsilon_{i}\big)>x\Big)\\
~&\notag=\E_{\sigma_{(k-1)h_n}}\bigg[\prod_{i=\lfloor k n h_n\rfloor+1}^{\lfloor (k+1)nh_n\rfloor}\P\Big(\epsilon_{i}>h_n^{1/2}\sigma_{(k-1)h_n}\big(x-h_n^{-1/2}(W_{t_i^n}-W_{kh_n})\big)|\mathcal{F}^X\Big)\bigg]\\
~&\notag=\E_{\sigma_{(k-1)h_n}}\bigg[\exp\Big(\sum_{i=\lfloor k n h_n\rfloor+1}^{\lfloor (k+1)nh_n\rfloor}\log\big(1-F_{\eta}\big(h_n^{1/2}\sigma_{(k-1)h_n}\big(x-h_n^{-1/2}(W_{t_i^n}-W_{kh_n})\big)\big)\big)\Big)\bigg]
\end{align}
by the tower rule for conditional expectations, and since $\epsilon_{i}\stackrel{iid}{\sim}F_{\eta}$. It holds that
\begin{align*}W_{t_i^n}-W_{kh_n}&=\sum_{j=1}^{i-\lfloor k n h_n\rfloor}\tilde U_j,~\tilde U_j\stackrel{iid}{\sim}\mathcal{N}(0,n^{-1}),j\ge 2,\tilde U_1\sim\mathcal{N}\big(0,t^n_{\lfloor k n h_n \rfloor +1}-kh_n\big)\,,\\
&U_j=h_n^{-1/2}\tilde U_j\stackrel{iid}{\sim}\mathcal{N}\big(0,(nh_n)^{-1}\big),j\ge 2,U_1\sim\mathcal{N}\big(0,h_n^{-1}\big(t^n_{\lfloor k n h_n \rfloor +1}-kh_n\big)\big)\,.\end{align*}
From \eqref{noise_dist}, and with a first-order Taylor expansion of $z\mapsto \log(1-z)$, we have that
\begin{align*}
\log\big(1-F_{\eta}(y)\big)\stackrel{\eqref{noise_dist}}{=}\log\big(1-\eta y(1+\KLEINO(1))\big)=-\eta y+\KLEINO(\eta y)=-\eta y_++\KLEINO(\eta y)\,,
\end{align*}
as $y\to 0$, where we add the positive part in the last equality, since $F_{\eta}(y)=0$, for any $y\le 0$. We obtain that
\begin{align*}
&\P_{\sigma_{(k-1)h_n}}\Big(h_n^{-1/2}\min_{i\in\mathcal{I}_{k}^n} \big(\epsilon_{i}+\sigma_{(k-1)h_n}(W_{t_i^n}-W_{kh_n})\big)>x\sigma_{(k-1)h_n}\Big)=\\
~&=\E_{\sigma_{(k-1)h_n}}\Big[\exp\Big(-h_n^{1/2}\sigma_{(k-1)h_n}\eta\sum_{i=\lfloor k n h_n \rfloor+1}^{\lfloor (k+1)nh_n\rfloor}\Big(x-\sum_{j=1}^{i-\lfloor k n h_n\rfloor}U_j\Big)_{+}(1+\KLEINO(1))\Big)\Big]\\
~&=\E_{\sigma_{(k-1)h_n}}\Big[\exp\Big(-h_n^{1/2}nh_n\sigma_{(k-1)h_n}\eta\int_0^1(B_t-x)_{-}\,dt\,(1+\KLEINO(1))\Big)\Big]\,.
\end{align*}
In the last equality, we use that the Riemann sums tend almost surely to the integral with a standard Brownian motion $(B_t)_{t\ge 0}$ in the integrand. Since the expression in the expectation is bounded, as a product of conditional probabilities by 1, we conclude with dominated convergence.
If $nh_n^{3/2}\to \infty$, we deduce that
\begin{align}
&\notag\P_{\sigma_{(k-1)h_n}}\Big(h_n^{-1/2}\min_{i\in\mathcal{I}_{k}^n} \big(\epsilon_{i}+\sigma_{(k-1)h_n}(W_{t_i^n}-W_{kh_n})\big)>x\sigma_{(k-1)h_n}\Big)=\P\big(\inf_{0\le t\le 1} B_t\ge x\big)\\
~&\notag+\E_{\sigma_{(k-1)h_n}}\Big[\1\big(\inf_{0\le t\le 1} B_t< x\big)\exp\Big(-h_n^{3/2}n\sigma_{(k-1)h_n}\eta\int_0^1(B_t-x)_{-}\,dt\,(1+\KLEINO(1))\Big)\Big]\\
~&=\P\big(\inf_{0\le t\le 1} B_t\ge x\big)+\P\big(\inf_{0\le t\le 1} B_t< x\big)\cdot \KLEINO(1)\,.\label{crucial}
\end{align}
We do not have a lower bound for $\int_0^1(B_t-x)_{-}\,dt$. However, using that the first entry time $T_x$ of $(B_t)$ in $x$, conditional on $\{\inf_{0\le t\le 1} B_t< x \}$, has a continuous conditional density $f(t|T_x<1)$, by Lemma \ref{lem3} and properties of the Brownian motion we obtain for any $\delta>0$ that
\begin{align*}
&\E_{\sigma_{(k-1)h_n}}\Big[\1\big(\inf_{0\le t\le 1} B_t< x\big)\exp\Big(-h_n^{3/2}n\sigma_{(k-1)h_n}\eta\int_0^1(B_t-x)_{-}\,dt\Big)\Big]\\
&\le \exp\big(-\big(h_n^{3/2}n\big)^{\delta}\sigma_{(k-1)h_n}\eta\big)\P(\inf_{0\le t\le 1} B_t< x)+\P\Big(\inf_{0\le t\le 1} B_t< x,\int_0^1(B_t-x)_{-}\,dt\le \big(h_n^{3/2}n\big)^{-1+\delta}\Big)\\
&\le \bigg(\hspace*{-.1cm}\exp\big(-\big(h_n^{3/2}n\big)^{\delta}\sigma_{(k-1)h_n}\eta\big)\hspace*{-.05cm}+\hspace*{-.05cm}\int_0^1 \P\Big(\int_s^1(B_t)_{-}\,dt\le \big(h_n^{3/2}n\big)^{-1+\delta}\Big) f(s|T_x<1)\,ds\hspace*{-.1cm}\bigg)\P(\inf_{0\le t\le 1} B_t< x)\\
&\le \bigg(\exp\big(-\big(h_n^{3/2}n\big)^{\delta}\sigma_{(k-1)h_n}\eta\big)+\int_0^1 \P\Big((1-s)\int_0^1(B_t)_{-}\,dt\le \big(h_n^{3/2}n\big)^{-1+\delta}\Big)\\
&\hspace*{7cm}\times f(s|T_x<1)\,ds\bigg)\,\P(\inf_{0\le t\le 1} B_t< x)\,.
\end{align*}
We focus on the second addend of the first factor, since the exponential term decays faster. It is bounded by a constant times
\begin{align*}
&\int_0^1 \P\Big((1-s)\int_0^1(B_t)_{-}\,dt\le \big(h_n^{3/2}n\big)^{-1+\delta}\Big)\,ds
\\ &\le \int_0^{1-b_n} \P\Big((1-s)\int_0^1(B_t)_{-}\,dt\le \big(h_n^{3/2}n\big)^{-1+\delta}\Big)\,ds+\int_{1-b_n}^1 \,ds\\
&\le \P\Big(b_n\int_0^1(B_t)_{-}\,dt\le \big(h_n^{3/2}n\big)^{-1+\delta}\Big)+b_n=\mathcal{O}\Big(\big(h_n^{3/2}nb_n^{-1}\big)^{-\frac{1+\delta}{3}}+b_n\Big)\,,
\end{align*}
for any sequence $(b_n)$, $b_n\in(0,1)$, where we use Lemma \ref{lem3}. Choosing $b_n$ which minimizes the order yields that 
\begin{align*}
\E_{\sigma_{(k-1)h_n}}\Big[\1\big(\inf_{0\le t\le 1} B_t< x\big)\exp\Big(-h_n^{3/2}n\sigma_{(k-1)h_n}\eta\int_0^1(B_t-x)_{-}\,dt\Big)\Big]=\P(\inf_{0\le t\le 1} B_t< x)\cdot R_n \,,
\end{align*}
almost surely, with a remainder which satisfies
\begin{align*}R_n=\mathcal{O}\Big(\big(h_n^{3/2}n\big)^{-\frac{1+\delta}{4}}\Big)\,.
\end{align*}
From the unconditional Lévy distribution of $T_x$, $f(s|T_x<1)$ is explicit, but we omit its precise form which does not influence the asymptotic order.
Under the condition $nh_n^{3/2}\to \infty$, the minimum of the Brownian motion over the interval hence dominates the noise in the distribution of local minima, different than for a choice $h_n\propto n^{-2/3}$. By the reflection principle, it holds that
\begin{align}\label{reflect}\P\big(-\inf_{0\le t\le 1} B_t\ge x\big)=\P\big(\sup_{0\le t\le 1} B_t\ge x\big)=2\P\big(B_1\ge x\big)=\P\big(|B_1|\ge x\big)\,,\end{align}
for $x\ge 0$.

Using the illustration of moments by integrals over tail probabilities we exploit this, and a completely analogous estimate for $\tilde m^*_{k-1,n}$, to approximate conditional expectations. This yields for all $k\in\{1,\ldots,h_n^{-1}-1\}$ that
\begin{align*}&\E_{\sigma_{(k-1)h_n}}\big[h_n^{-1/2}\tilde m_{k,n}\big]=\\
&=\int_0^{\infty}\P_{\sigma_{(k-1)h_n}}\big(h_n^{-1/2}\tilde m_{k,n}>x\big)\,dx-\int_0^{\infty}\P_{\sigma_{(k-1)h_n}}\big(-h_n^{-1/2}\tilde m_{k,n}>x\big)\,dx\\
&=-\int_0^{\infty}\P_{\sigma_{(k-1)h_n}}\big(\sigma_{(k-1)h_n}\sup_{0\le t\le 1} B_t>x\big)\,dx+\KLEINO_{\P}(1)\\
&=-\int_0^{\infty}\P_{\sigma_{(k-1)h_n}}\big(\sigma_{(k-1)h_n}|B_1|>x\big)\,dx+\KLEINO_{\P}(1)\\
&=-\E_{\sigma_{(k-1)h_n}}\big[\sigma_{(k-1)h_n}| B_1|\big]+\KLEINO_{\P}(1)\\
&=-\sqrt{\frac{2}{\pi}}\sigma_{(k-1)h_n}+\KLEINO_{\P}(1)\,.
\end{align*}
We used \eqref{reflect}. An analogous computation yields the same result for $\tilde m^*_{k-1,n}$:
\begin{align*}\E_{\sigma_{(k-1)h_n}}\big[h_n^{-1/2}\tilde m^*_{k-1,n}\big]=-\sqrt{\frac{2}{\pi}}\sigma_{(k-1)h_n}+\KLEINO_{\P}(1)\,.
\end{align*}
For the second conditional moments, we obtain for all $k\in\{1,\ldots,h_n^{-1}-1\}$ that
\begin{align*}\E_{\sigma_{(k-1)h_n}}\big[h_n^{-1}\big(\tilde m_{k,n}\big)^2\big]&=2\int_0^{\infty}x\,\P_{\sigma_{(k-1)h_n}}\big(|h_n^{-1/2}\tilde m_{k,n}|>x\big)\,dx\\
&=2\int_0^{\infty}x\,\P_{\sigma_{(k-1)h_n}}\big(\sigma_{(k-1)h_n}\sup_{0\le t\le 1} B_t>x\big)\,dx+\KLEINO_{\P}(1)\\
&=2\int_0^{\infty}x\,\P_{\sigma_{(k-1)h_n}}\big(\sigma_{(k-1)h_n}|B_1|>x\big)\,dx+\KLEINO_{\P}(1)\\
&=\sigma_{(k-1)h_n}^2+\KLEINO_{\P}(1)\,.
\end{align*}
The last identity uses the illustration of the second moment of the normal distribution as an integral over tail probabilities. An analogous computation yields that
\begin{align*}\E_{\sigma_{(k-1)h_n}}\big[h_n^{-1}\big(\tilde m^*_{k-1,n}\big)^2\big]=\sigma_{(k-1)h_n}^2+\KLEINO_{\P}(1)\,.
\end{align*}
Inserting the identities for the conditional moments in \eqref{tower} yields that
\begin{align*}\E\big[h_n^{-1}\big(\tilde m_{k,n}-\tilde m^*_{k-1,n})^2\big]&=2\Big(1-\frac{2}{\pi}\Big)\E[\sigma_{(k-1)h_n}^2]+\KLEINO(1)\,,
\end{align*}
such that
\begin{align*}
\E\big[\hat\sigma^2_{\tau-}-\sigma_{\tau-}^2\big]&=\frac{\pi}{2(\pi-2)K_n}\sum_{k=(\lfloor h_n^{-1}\tau\rfloor-K_n)\vee 1}^{\lfloor h_n^{-1}\tau\rfloor-1}h_n^{-1}\E\Big[\big(\tilde m_{k,n}-\tilde m^*_{k-1,n})^2\Big]-\E[\sigma_{\tau-}^2]+\KLEINO(1)\\
&=\frac{1}{K_n}\sum_{k=(\lfloor h_n^{-1}\tau\rfloor-K_n)\vee 1}^{\lfloor h_n^{-1}\tau\rfloor-1}\E[\sigma_{(k-1)h_n}^2-\sigma_{\tau-}^2]+\KLEINO(1)=\KLEINO(1)\,.
\end{align*}
This proves \eqref{expsimple}. Since the next step shows that the variance of the estimator tends to zero, consistency holds true.\\[.2cm] 
{\bf{Step 4}}\\[.1cm]
We determine the asymptotic variance of the estimator. Illustrating moments as integrals over tail probabilities, with the analogous approximation as above, we obtain for all $k\in\{1,\ldots,h_n^{-1}-1\}$ that
\begin{align*}&\var_{\sigma_{(k-1)h_n}}\big(\tilde m_{k,n}^2\big)=\E_{\sigma_{(k-1)h_n}}\big[\tilde m_{k,n}^4\big]-\Big(\E_{\sigma_{(k-1)h_n}}\big[\tilde m_{k,n}^2\big]\Big)^2=2\sigma_{(k-1)h_n}^4h_n^2+\KLEINO_{\P}(h_n^2),\\
&\cov_{\sigma_{(k-1)h_n}}\big(\tilde m_{k,n}^2,\tilde m_{k,n}\tilde m^*_{k-1,n}\big)=\E_{\sigma_{(k-1)h_n}}\big[\tilde m_{k,n}^3\big]\E_{\sigma_{(k-1)h_n}}\big[\tilde m^*_{k-1,n}\big]\\&\hspace*{5.25cm}-\E_{\sigma_{(k-1)h_n}}\big[\tilde m_{k,n}^2\big]\E_{\sigma_{(k-1)h_n}}\big[\tilde m_{k,n}\big]\E_{\sigma_{(k-1)h_n}}\big[\tilde m^*_{k-1,n}\big]\\
&\hspace*{4.85cm}=\frac{2}{\pi}\,\sigma_{(k-1)h_n}^4h_n^2+\KLEINO_{\P}(h_n^2),\\
&\var_{\sigma_{(k-1)h_n}}\Big(\tilde m_{k,n}\tilde m^*_{k-1,n}\Big)=\E_{\sigma_{(k-1)h_n}}\big[\tilde m_{k,n}^2\big]\E_{\sigma_{(k-1)h_n}}\big[\big(\tilde m^*_{k-1,n}\big)^2\big]\\
&\hspace*{4.35cm}-\big(\E_{\sigma_{(k-1)h_n}}\big[\tilde m_{k,n}\big]\E_{\sigma_{(k-1)h_n}}\big[\tilde m^*_{k-1,n}\big]\big)^2\\
&\hspace*{3.95cm}=\sigma_{(k-1)h_n}^4\Big(1-\frac{4}{\pi^2}\Big)h_n^2+\KLEINO_{\P}(h_n^2).
\end{align*}
We have used the first four moments of the half-normal distribution and their illustration via integrals over tail probabilities. The dependence structure between  $\tilde m_{k,n}$ and $\tilde m^*_{k,n}$ also affects the variance of $\hat\sigma^2_{\tau-}$. We perform approximation steps for covariances similar as for the moments of local minima above, using that
\begin{align*}&h_n^{-1}\cov_{\sigma_{(k-1)h_n}}\big(\tilde m_{k,n},\tilde m^*_{k,n}\big)=\int_{-\infty}^{\infty}\int_{-\infty}^{\infty}\Big(\P_{\sigma_{(k-1)h_n}}\big(h_n^{-1/2}\tilde m_{k,n}>x, h_n^{-1/2}\tilde m^*_{k,n}>y\big)\\
&\hspace*{5.5cm}-\P_{\sigma_{(k-1)h_n}}\big(h_n^{-1/2}\tilde m_{k,n}>x\big)\P_{\sigma_{(k-1)h_n}}\big(h_n^{-1/2}\tilde m^*_{k,n}>y\big)\Big)\,dxdy\\
&=\int_{0}^{\infty}\int_{0}^{\infty}\Big(\P_{\sigma_{(k-1)h_n}}\big(\sigma_{(k-1)h_n}\sup_{0\le t\le 1}B_t>x, \sigma_{(k-1)h_n}\big(\sup_{0\le t\le 1}B_t-B_1\big)>y\big)\\
&\hspace*{1cm}-\P_{\sigma_{(k-1)h_n}}\big(\sigma_{(k-1)h_n}\sup_{0\le t\le 1}B_t>x\big)\P_{\sigma_{(k-1)h_n}}\big(\sigma_{(k-1)h_n}\hspace*{-.05cm}\big(\sup_{0\le t\le 1}B_t-B_1\big)>y\big)\Big)dxdy\hspace*{-.05cm}+\hspace*{-.05cm}\KLEINO_{\P}(1).
\end{align*}
This shows that the joint distribution of $(\tilde m_{k,n},\tilde m^*_{k,n})$ relates to the distribution of the minimum and the difference between minimum and endpoint of Brownian motion over an interval, or equivalently the distribution of the maximum and the difference between maximum and endpoint. The latter is readily obtained from the joint density of maximum and endpoint which is a well-known result on stochastic processes, see e.g.\ \cite{shepp}. Utilizing this, we obtain for all $k\in\{1,\ldots,h_n^{-1}-1\}$ that
\[\cov_{\sigma_{(k-1)h_n}}\big(\tilde m_{k,n}\,,\,\tilde m^*_{k,n}\big)=\Big(\frac{1}{2}-\frac{2}{\pi}\Big)h_n \,\sigma_{(k-1)h_n}^2(1+\mathcal{O}_{\P}(h_n^{\alpha}))+\KLEINO_{\P}\big(h_n\big)\,.\]
The additional remainder of order $h_n^{\alpha}$ in probability is due to the different approximations of $(\sigma_t)$ in $\tilde m_{k,n}$ and $\tilde m^*_{k,n}$. This implies for all $k\in\{1,\ldots,h_n^{-1}-1\}$ that
\begin{align*}
&\cov_{\sigma_{(k-1)h_n}}\big(\tilde m_{k,n}\tilde m^*_{k-1,n},\tilde m_{k+1,n}\tilde m^*_{k,n}\big)=\\
&~\Big(\E_{\sigma_{(k-1)h_n}}\big[\tilde m_{k,n}\tilde m^*_{k,n}\big]-\E_{\sigma_{(k-1)h_n}}\big[\tilde m_{k,n}\big]\E_{\sigma_{(k-1)h_n}}\big[\tilde m^*_{k,n}\big]\Big) \E_{\sigma_{(k-1)h_n}}\big[\tilde m^*_{k-1,n}\big]\E\big[\tilde m_{k+1,n}\big]\\
&=\sigma_{(k-1)h_n}^4\Big(\frac{1}{\pi}-\frac{4}{\pi^2}\Big)h_n^2+\KLEINO_{\P}\big(h_n^2\big)\,.
\end{align*}
With analogous steps, we deduce two more covariances which contribute to the asymptotic variance:
\begin{align*}
&\cov_{\sigma_{(k-1)h_n}}\big(\tilde m_{k,n}^2,\big(\tilde m^*_{k,n}\big)^2\big)=-h_n^2\frac{\sigma_{(k-1)h_n}^4}{2}+\KLEINO_{\P}\big(h_n^2\big)\,,\\
&\cov_{\sigma_{(k-1)h_n}}\big(\big(\tilde m^*_{k,n}\big)^2,m_k\tilde m^*_{k-1,n}\big)=-h_n^2\frac{2}{3\pi}\sigma_{(k-1)h_n}^4+\KLEINO_{\P}\big(h_n^2\big)\,.
\end{align*}
All covariance terms which enter the asymptotic variance are of one of these forms. For the conditional variance given $\sigma^2_{\tau-}$, we obtain that
\begin{align*}&\var_{\sigma^2_{\tau-}}\negthinspace\negthinspace\big(\hat\sigma^2_{\tau-}\big)=\frac{1}{K_n^2}\frac{\pi^2}{4(\pi-2)^2}\bigg(\sum_{k=(\lfloor h_n^{-1}\tau\rfloor-K_n)\vee 1}^{\lfloor h_n^{-1}\tau\rfloor-1}h_n^{-2}\var_{\sigma^2_{\tau-}}\negthinspace\negthinspace\Big(\tilde m_{k,n}^2+(\tilde m^*_{k,n})^2-2\tilde m_{k,n}\tilde m^*_{k-1,n}\Big)\\
&-\negthinspace\sum_{k=(\lfloor h_n^{-1}\tau\rfloor-K_n)\vee 2}^{\lfloor h_n^{-1}\tau\rfloor-1}\negthinspace\negthinspace 4h_n^{-2}\cov_{\sigma^2_{\tau-}}\negthinspace\negthinspace\Big(\tilde m_{k,n}\tilde m^*_{k-1,n}\,,\,\tilde m_{k-1,n}^2\negthinspace+\negthinspace(\tilde m^*_{k-1,n})^2\negthinspace-\negthinspace2\tilde m_{k-1,n}\tilde m^*_{k-2,n}\Big)\negthinspace\bigg)\negthinspace+\negthinspace\KLEINO_{\P}\big(K_n^{-1}\big)\\
&=\frac{1}{K_n^2}\frac{\pi^2}{4(\pi-2)^2}\bigg(\sum_{k=(\lfloor h_n^{-1}\tau\rfloor-K_n)\vee 1}^{\lfloor h_n^{-1}\tau\rfloor-1}h_n^{-2}\Big(2\var_{\sigma^2_{\tau-}}\negthinspace\negthinspace\big(\tilde m_{k,n}^2\big)+4\var_{\sigma^2_{\tau-}}\negthinspace\big(\tilde m_{k,n}\tilde m^*_{k-1,n}\big)\\
&\hspace*{.0cm}+2\cov_{\sigma^2_{\tau-}}\negthinspace\negthinspace\big(\tilde m_{k,n}^2,(\tilde m^*_{k,n})^2\big)-4\cov_{\sigma^2_{\tau-}}\negthinspace\negthinspace\big(\tilde m_{k,n}^2, \tilde m_{k,n}\tilde m^*_{k-1,n}\big)-4\cov_{\sigma^2_{\tau-}}\negthinspace\negthinspace\big((\tilde m^*_{k,n})^2, \tilde m_{k,n}\tilde m^*_{k-1,n}\big)\Big)\\
&\hspace*{.0cm}+\negthinspace\sum_{k=(\lfloor h_n^{-1}\tau\rfloor-K_n)\vee 2}^{\lfloor h_n^{-1}\tau\rfloor-1}\negthinspace\negthinspace 4h_n^{-2}\Big(2\cov_{\sigma^2_{\tau-}}\negthinspace\negthinspace\big(\tilde m_{k,n}\tilde m^*_{k-1,n}, \tilde m_{k-1,n}\tilde m^*_{k-2,n}\big)\negthinspace-\negthinspace\cov_{\sigma^2_{\tau-}}\negthinspace\negthinspace\big(\tilde m_{k,n}\tilde m^*_{k-1,n},\tilde m_{k-1,n}^2\big)\\
&\hspace*{4.75cm}-\cov_{\sigma^2_{\tau-}}\negthinspace\negthinspace\big(\tilde m_{k,n}\tilde m^*_{k-1,n},(\tilde m^*_{k-1,n})^2\big)\Big)
\bigg)+\KLEINO_{\P}\big(K_n^{-1}\big)\\
&=\frac{1}{K_n}\frac{\pi^2}{4(\pi-2)^2}\sigma_{\tau-}^4\Big(8-\frac{16}{\pi^2}-1-\frac{8}{\pi}+\frac{8}{3\pi}+2\Big(\frac{4}{3\pi}-\frac{16}{\pi^2}\Big)\Big)+\KLEINO_{\P}\big(K_n^{-1}\big)\\
&=\frac{1}{K_n}\frac{1}{(\pi-2)^2}\Big(\frac{7\pi^2}{4}-\frac{2\pi}{3}-12\Big)\sigma^4_{\tau-}+\KLEINO_{\P}\big(K_n^{-1}\big)\,.
\end{align*}
{\bf{Step 5}}\\[.1cm]
For a central limit theorem, the squared bias needs to be asymptotically negligible compared to the variance, which is satisfied for $K_n=\KLEINO(h_n^{-2\alpha/(1+2\alpha)})$. By the existence of higher moments of $\tilde m_{k,n}$ and $\tilde m^*_{k-1,n}$, a Lyapunov-type condition is straightforward, such that asymptotic normality conditional on $\sigma^2_{\tau-}$ is implied by a classical central limit theorem for m-dependent triangular arrays as the one by \cite{berk}. A feasible central limit theorem is implied by this conditional asymptotic normality in combination with  $\mathcal{F}^X$-stable convergence. For the stability, we show that $\alpha_n=K_n^{1/2}\big(\hat\sigma_{\tau-}^2-\sigma_{\tau-}^2\big)$ satisfy
\begin{align}\label{stable}\E\left[Z g(\alpha_n)\right]\rightarrow \E\left[Z g(\alpha)\right]=\E[Z]\E\left[g(\alpha)\right]\,,\end{align}
for any $\mathcal{F}^X$-measurable bounded random variable $Z$ and continuous bounded function $g$, where
\begin{align}\alpha=\sigma_{\tau-}^2\frac{1}{(\pi-2)}\sqrt{\frac{7\pi^2}{4}-\frac{2\pi}{3}-12}~ U \,,\end{align}
with $U$ a standard normally distributed random variable which is independent of $\mathcal{F}^X$. By the above approximations it suffices to prove this for the statistics based on $\tilde m_{k,n}$ and $\tilde m^*_{k-1,n}$ from \eqref{eqlem1}, and $Z$ measurable w.r.t.\ $\sigma(\int_0^t\sigma_s\,dW_s,0\le t\le 1)$. Set
\begin{align*}A_n=[\tau-(K_n+1)h_n,\tau]\,,\,\tilde X(n)_t=\int_0^t\1_{A_n}(s)\sigma_{\lfloor sh_n^{-1}\rfloor h_n}\,dW_s\,,
\bar X(n)_t=X_t-\tilde X(n)_t\,.
\end{align*}
Denote with $\mathcal{H}_n$ the $\sigma$-field generated by $\bar X(n)_t$ and $\mathcal{F}^X_0$. The sequence $\big(\mathcal{H}_n\big)_{n\in\N}$ is isotonic with limit $\bigvee_n \mathcal{H}_n=\sigma(\int_0^t\sigma_s\,dW_s,0\le t\le 1)$. Since $\E[Z|\mathcal{H}_n]\rightarrow Z$ in $L^1(\P)$ as $n\rightarrow\infty$, it is enough to show that $\E[Zg(\alpha_n)]\rightarrow \E[Z]\E[g(\alpha)]$, for $Z$ being $\mathcal{H}_{n_0}$-measurable for some $n_0\in\N$. Observe that $\alpha_n$ includes only increments of local minima based on $\tilde X(n)_t$, which are uncorrelated from those of $\bar X(n)_t$. For all $n\ge n_0$, we hence obtain that $\E[Zg(\alpha_n)]=\E[Z]\E[g(\alpha_n)]\rightarrow \E[Z]\E[g(\alpha)]$ by a standard central limit theorem. This shows \eqref{stable} and completes the proof of \eqref{spotclt}.
\subsubsection{Proof of Proposition \ref{propvola2}}
For the quarticity estimator \eqref{quartestimator}, when $\lfloor h_n^{-1}\tau\rfloor>K_n$, we have that
\begin{align*}
\E\big[\widehat{{\sigma^4_{\tau}}}_--\sigma^4_{\tau-}\big]&=\frac{\pi}{4(3\pi-8)K_n}\sum_{k=(\lfloor h_n^{-1}\tau\rfloor-K_n)\vee 1}^{\lfloor h_n^{-1}\tau\rfloor-1}h_n^{-2}\E\Big[\tilde m_{k,n}^4+(\tilde m^*_{k-1,n})^4-4 \tilde m_{k,n}^3\tilde m^*_{k-1,n}\\
&\hspace*{1.65cm}-4 \tilde m_{k,n}(\tilde m^*_{k-1,n})^3+6 \tilde m_{k,n}^2(\tilde m^*_{k-1,n})^2\Big]-\E[\sigma^4_{\tau-}]+\mathcal{O}\big(h_n^{\alpha\wedge 1/2}\big)\\
&=\Big(\frac{\pi}{4(3\pi-8)}\big(6-16/\pi-16/\pi+6\big)-1\Big)\E[\sigma^4_{\tau-}]+\KLEINO(1)\\
&=\KLEINO(1)\,,
\end{align*}
by using the same moments as in the computation of the asymptotic variance. We can bound its variance by
\begin{align*}
\var\big(\widehat{{\sigma^4_{\tau}}}_-\big)&\le\frac{\pi^2}{16(3\pi-8)^2K_n^2}2K_nh_n^{-4}\,\var\Big(\big(\tilde m_{k,n}-\tilde m^*_{k-1,n}\big)^4\Big)+\KLEINO\big(K_n^{-1}\big)\\
&\le \frac{1}{K_n}\frac{\pi^2}{8(3\pi-8)^2}h_n^{-4}\,\E\Big[\big(\tilde m_{k,n}-\tilde m^*_{k-1,n}\big)^8\Big]+\KLEINO\big(K_n^{-1}\big)\\
&\le \frac{1}{K_n}\frac{\pi^2}{8(3\pi-8)^2}h_n^{-4}\,256\,\E\Big[\tilde m_{k,n}^8\Big]+\KLEINO\big(K_n^{-1}\big)=\mathcal{O}(K_n^{-1})\,,\\
\end{align*}
what readily implies Proposition \ref{propvola2}.

\subsection{Asymptotics of the truncated spot volatility estimation with jumps}
Denote by 
\begin{align*}D^X_k:=m_{k,n}-m_{k-1,n}, ~k=1,\ldots,h_n^{-1}-1\,,
\end{align*}
the differences of local minima based on the observations \eqref{lomn}, with the general semimartingale \eqref{smdecomp} with jumps.
Denote by 
\begin{align*}D^C_k:=\tilde m_{k,n}-\tilde m_{k-1,n}^*, ~k=1,\ldots,h_n^{-1}-1\,,
\end{align*}
the differences of the unobservable local minima considered in Section \ref{sec:volacont}. In particular, the statistics $D^C_k$ are based only on the continuous part $(C_t)$ in \eqref{smdecomp} such that the jumps are eliminated. Theorem \ref{propvola3} is implied by Proposition \ref{propvola2}, if we can show that
\begin{align*}\frac{\pi}{2(\pi-2)K_n}\sum_{k=(\lfloor h_n^{-1}\tau\rfloor-K_n)\vee 1}^{\lfloor h_n^{-1}\tau\rfloor-1}h_n^{-1}\Big(\big(D^X_k\big)^2\1_{\{|D^X_k|\le u_n\}}-\big(D^C_k\big)^2\Big)=\mathcal{O}_{\P}\big(h_n^{\frac{\alpha}{2\alpha+1}}\big)=\KLEINO_{\P}\big(K_n^{-1/2}\big)\,.\
\end{align*}
We decompose this difference of the truncated estimator, which is based on the available observations with jumps, and the non-truncated estimator, which uses non-available observations without jumps, in the following way:
\begin{align*}&\frac{\pi}{2(\pi-2)K_n}\sum_{k=(\lfloor h_n^{-1}\tau\rfloor-K_n)\vee 1}^{\lfloor h_n^{-1}\tau\rfloor-1}h_n^{-1}\Big(\big(D^X_k\big)^2\1_{\{|D^X_k|\le u_n\}}-\big(D^C_k\big)^2\Big)\\
&=\frac{\pi}{2(\pi-2)K_n}\sum_{k=(\lfloor h_n^{-1}\tau\rfloor-K_n)\vee 1}^{\lfloor h_n^{-1}\tau\rfloor-1}h_n^{-1}\bigg(\1_{\{|D^C_k|> c \hspace*{.05em} u_n\}}\Big(\big(D^X_k\big)^2\1_{\{|D^X_k|\le u_n\}}-\big(D^C_k\big)^2\Big)\\
&\hspace*{6cm} +\1_{\{|D^C_k|\le  c \hspace*{.05em} u_n\}}\1_{\{|D^X_k|\le u_n\}}\Big(\big(D^X_k\big)^2-\big(D^C_k\big)^2\Big)\\
&\hspace*{6cm} -\1_{\{|D^C_k|\le  c \hspace*{.05em} u_n\}}\1_{\{|D^X_k|> u_n\}}\big(D^C_k\big)^2\bigg)\,,
\end{align*}
with some arbitrary constant $c\in(0,1)$. Without loss of generality we can set $\beta=1$ in this proof, i.e.\ $u_n=h_n^{\kappa}$. We consider the three addends which are different error terms by
\begin{enumerate}[itemsep=2pt,parsep=2pt]
\item large absolute statistics based on the continuous part  $(C_t)$;
\item non-truncated statistics which contain (small) jumps;
\item the truncation of also the continuous parts in statistics $(D_k^X)$ which exceed the threshold;
\end{enumerate}
separately. The probability $\P(|D^C_k|> c  \hspace*{.05em} u_n)$ can be bounded using the estimate from \eqref{crucial} and Gaussian tail bounds. Observe that the remainder in \eqref{crucial} is non-negative. This yields that for some $y>0$, we have that 
\[\P\Big(h_n^{-1/2}\big|\tilde m_{k,n}\big|>y\Big)\le \P\Big(\sup_{0\le t\le 1} B_t >y\Big)\,,\]
what is intuitive, since the errors $(\epsilon_i)$ are non-negative. We apply the triangular inequality and then Hölder's inequality to the expectation of the absolute first error term and obtain for any $p\in \N$ that
\begin{align*}&\frac{\pi}{2(\pi-2)K_n}\E\bigg[\bigg|\sum_{k=(\lfloor h_n^{-1}\tau\rfloor-K_n)\vee 1}^{\lfloor h_n^{-1}\tau\rfloor-1}h_n^{-1}\1_{\{|D^C_k|> c \hspace*{.05em} u_n\}}\Big(\big(D^X_k\big)^2\1_{\{|D^X_k|\le u_n\}}-\big(D^C_k\big)^2\Big)\bigg|\bigg]\\
&\le\frac{\pi}{2(\pi-2)K_n}\sum_{k=(\lfloor h_n^{-1}\tau\rfloor-K_n)\vee 1}^{\lfloor h_n^{-1}\tau\rfloor-1}h_n^{-1}\E\Big[\1_{\{|D^C_k|> c \hspace*{.05em}u_n\}}\Big|\big(D^X_k\big)^2\1_{\{|D^X_k|\le u_n\}}-\big(D^C_k\big)^2\Big|\Big]\\
&\le \frac{\pi}{2(\pi-2)K_n}\sum_{k=(\lfloor h_n^{-1}\tau\rfloor-K_n)\vee 1}^{\lfloor h_n^{-1}\tau\rfloor-1}h_n^{-1}\Big(\P\big(|D^C_k|> c \hspace*{.05em} u_n\big)2\,\big(u_n^4+\E\big[\big(D^C_k\big)^4\big]\big)\Big)^{1/2}\\
&\le \frac{\pi}{2(\pi-2)K_n}\sum_{k=(\lfloor h_n^{-1}\tau\rfloor-K_n)\vee 1}^{\lfloor h_n^{-1}\tau\rfloor-1}h_n^{-1}\Big(\P\big(h_n^{-1/2}|D^C_k|> c h_n^{\kappa-1/2}\big)\Big)^{1/2}\sqrt{2}\, u_n^2\\
&\le \frac{\pi}{2(\pi-2)K_n}\sum_{k=(\lfloor h_n^{-1}\tau\rfloor-K_n)\vee 1}^{\lfloor h_n^{-1}\tau\rfloor-1}h_n^{-1}\Big(2\,\P\Big(|B_1|> \frac{c}{2} h_n^{\kappa-1/2}\Big)\Big)^{1/2}\sqrt{2}\, u_n^2\\
&\le \frac{\sqrt{2}\pi}{(\pi-2)K_n}\sum_{k=(\lfloor h_n^{-1}\tau\rfloor-K_n)\vee 1}^{\lfloor h_n^{-1}\tau\rfloor-1}h_n^{2\kappa-1}\exp\Big(- \frac{c^2}{4} h_n^{2\kappa-1}\Big)\\
&=\mathcal{O}\Big(h_n^{(-p+1)(2\kappa-1)}\Big)=\KLEINO\big(h_n^{\frac{\alpha}{2\alpha+1}}\big)\,.
\end{align*}
Since $2\kappa-1<0$ and $p$ arbitrarily large, we conclude that the first error term is asymptotically negligible. We will use the elementary inequalities
\begin{align*}D^X_k&=\min_{i\in\mathcal{I}_k^n}\big(C_{\frac{i}{n}}+J_{\frac{i}{n}}+\epsilon_i\big)-\min_{i\in\mathcal{I}_{k-1}^n}\big(C_{\frac{i}{n}}+J_{\frac{i}{n}}+\epsilon_i\big)\\
&\le \min_{i\in\mathcal{I}_k^n}\big(C_{\frac{i}{n}}+\epsilon_i\big)+\max_{i\in\mathcal{I}_k^n}J_{\frac{i}{n}}-\min_{i\in\mathcal{I}_{k-1}^n}\big(C_{\frac{i}{n}}+\epsilon_i\big)-\min_{i\in\mathcal{I}_{k-1}^n}J_{\frac{i}{n}}\\
&=D^C_k+ \max_{i\in\mathcal{I}_k^n}J_{\frac{i}{n}}-\min_{i\in\mathcal{I}_{k-1}^n}J_{\frac{i}{n}}+\mathcal{O}_{\P}\big(h_n^{\alpha\wedge 1/2}\big)\,,
\end{align*}
and
\begin{align*}D^X_k&=\min_{i\in\mathcal{I}_k^n}\big(C_{\frac{i}{n}}+J_{\frac{i}{n}}+\epsilon_i\big)-\min_{i\in\mathcal{I}_{k-1}^n}\big(C_{\frac{i}{n}}+J_{\frac{i}{n}}+\epsilon_i\big)\\
&\ge  \min_{i\in\mathcal{I}_k^n}\big(C_{\frac{i}{n}}+\epsilon_i\big)+\min_{i\in\mathcal{I}_k^n}J_{\frac{i}{n}}-\min_{i\in\mathcal{I}_{k-1}^n}\big(C_{\frac{i}{n}}+\epsilon_i\big)-\max_{i\in\mathcal{I}_{k-1}^n}J_{\frac{i}{n}}\\
&=D^C_k+ \min_{i\in\mathcal{I}_k^n}J_{\frac{i}{n}}-\max_{i\in\mathcal{I}_{k-1}^n}J_{\frac{i}{n}}+\mathcal{O}_{\P}\big(h_n^{\alpha\wedge 1/2}\big)\,.
\end{align*}
Therefore, we can bound $|D^X_k-D^C_k|$ by
\begin{align*}
\sup_{\substack{i\in\mathcal{I}_{k}^n, j\in\mathcal{I}_{k-1}^n}}|J_{\frac{i}{n}}-J_{\frac{j}{n}}|&\le \sup_{\substack{s\in[kh_n,(k+1)h_n], t\in[(k-1)h_n,kh_n]}}|J_{s}-J_{t}|\\
&\le \sup_{\substack{s\in[kh_n,(k+1)h_n]}}|J_{s}-J_{kh_n}|+ \sup_{\substack{t\in[(k-1)h_n,kh_n]}}|J_{kh_n}-J_{t}|\,,
\end{align*}
and the remainder term of the approximation for the continuous part which is $\mathcal{O}_{\P}\big(h_n^{\alpha\wedge 1/2}\big)$. Since the compensated small jumps of a semimartingale admit a martingale structure, Doob's inequality for c\`{a}dl\`{a}g $L_2$-martingales can be used to bound these suprema. Based on these preliminaries, we obtain for the expected absolute value of the second error term that
\begin{align*}
&\frac{\pi}{2(\pi-2)K_n}\E\bigg[\bigg|\sum_{k=(\lfloor h_n^{-1}\tau\rfloor-K_n)\vee 1}^{\lfloor h_n^{-1}\tau\rfloor-1}h_n^{-1}\1_{\{|D^C_k|\le  c \hspace*{.05em} u_n\}}\1_{\{|D^X_k|\le u_n\}}\Big(\big(D^X_k\big)^2-\big(D^C_k\big)^2\Big)\bigg|\bigg]\\
&\le \frac{\pi}{2(\pi-2)K_n}\sum_{k=(\lfloor h_n^{-1}\tau\rfloor-K_n)\vee 1}^{\lfloor h_n^{-1}\tau\rfloor-1}h_n^{-1}\E\Big[\1_{\{|D^C_k|\le  c \hspace*{.05em} u_n\}}\1_{\{|D^X_k|\le u_n\}}\Big|\big(D^X_k\big)^2-\big(D^C_k\big)^2\Big|\Big]\\
&\lesssim \frac{1}{K_n}\sum_{k=(\lfloor h_n^{-1}\tau\rfloor-K_n)\vee 1}^{\lfloor h_n^{-1}\tau\rfloor-1}h_n^{-1}\E\Big[\sup_{\substack{i\in\mathcal{I}_{k}^n, j\in\mathcal{I}_{k-1}^n}}|J_{\frac{i}{n}}-J_{\frac{j}{n}}|^2\wedge (1+c)^2u_n^2\Big]\\
&\lesssim \frac{1}{K_n}\sum_{k=(\lfloor h_n^{-1}\tau\rfloor-K_n)\vee 1}^{\lfloor h_n^{-1}\tau\rfloor-1}h_n^{-1}\E\Big[\sup_{\substack{t\in [kh_n,(k+1)h_n]}}|J_{t}-J_{kh_n}|^2\wedge u_n^2\Big]\\
&\lesssim \frac{1}{K_n}\sum_{k=(\lfloor h_n^{-1}\tau\rfloor-K_n)\vee 1}^{\lfloor h_n^{-1}\tau\rfloor-1}h_n^{-1}\E\Big[|J_{(k+1)h_n}-J_{kh_n}|^2\wedge u_n^2\Big]\\
&=\mathcal{O}\big(u_n^{2-r}\big)\,.
\end{align*}
Applying the elementary inequalities from above, a cross term in the upper bound for $\big(D^X_k\big)^2-\big(D^C_k\big)^2$ is of smaller order and directly neglected. It can be handled using the Cauchy-Schwarz inequality. In the last step, we adopt a bound on the expected absolute thresholded jump increments from Equation (54) in \cite{aitjac10}. For the negligibility of the second error term, we thus get the condition that
\begin{align}\label{kappa1}\kappa(2-r)\ge \frac{\alpha}{1+2\alpha}\,.\end{align}
Doob's inequality yields as well that
\begin{align*}
\P\big(\sup_{\substack{t\in [kh_n,(k+1)h_n]}}|J_{t}-J_{kh_n}|\ge (1-c) u_n\big)&\le \frac{\E\big[\big|J_{(k+1)h_n}-J_{kh_n}\big|^{r\wedge 1}\big]}{\big((1-c)u_n\big)^{r\wedge 1}}+\mathcal{O}(h_n)\\
&=\mathcal{O}\big(h_n u_n^{-r}\big)\,.
\end{align*}
For this upper bound, we decomposed the jumps in the sum of large jumps and the martingale of compensated small jumps to which we apply Doob's inequality. We derive the following estimate for the expectation of the third (absolute) error term
\begin{align*}
&\frac{\pi}{2(\pi-2)K_n}\sum_{k=(\lfloor h_n^{-1}\tau\rfloor-K_n)\vee 1}^{\lfloor h_n^{-1}\tau\rfloor-1}h_n^{-1}\E\Big[\1_{\{|D^C_k|\le  c  \hspace*{.05em} u_n\}}\1_{\{|D^X_k|> u_n\}}\big(D^C_k\big)^2\Big]\\
&\le \frac{\pi}{2(\pi-2)K_n}\sum_{k=(\lfloor h_n^{-1}\tau\rfloor-K_n)\vee 1}^{\lfloor h_n^{-1}\tau\rfloor-1}h_n^{-1}\E\Big[\1_{\{2\sup_{\substack{s\in[(k-1)h_n,(k+1)h_n]}}|J_{s}-J_{kh_n}|\ge (1-c) u_n\}}\big(D^C_k\big)^2\Big]\\
&\lesssim \frac{1}{K_n}\sum_{k=(\lfloor h_n^{-1}\tau\rfloor-K_n)\vee 1}^{\lfloor h_n^{-1}\tau\rfloor-1}h_n^{-1}\P\Big(\sup_{\substack{t\in [kh_n,(k+1)h_n]}}|J_{t}-J_{kh_n}|\ge (1-c) u_n\Big)\E\Big[\big(D^C_k\big)^2\Big]\\
&\lesssim \frac{1}{K_n}\sum_{k=(\lfloor h_n^{-1}\tau\rfloor-K_n)\vee 1}^{\lfloor h_n^{-1}\tau\rfloor-1}\bigg(\frac{\E\big[\big|J_{(k+1)h_n}-J_{kh_n}\big|^{r\wedge 1}\big]}{\big((1-c)u_n\big)^{r\wedge 1}}+\mathcal{O}(h_n)\bigg)\\
&=\mathcal{O}\big(h_n u_n^{-r}\big)\,.
\end{align*}
For the negligibility of the third error term, we thus get the condition that
\begin{align}\label{kappa2}1-\kappa r\ge \frac{\alpha}{1+2\alpha}\,.\end{align}
Since under the conditions of Theorem \ref{propvola3}, \eqref{kappa1} and \eqref{kappa2} are satisfied, the proof is finished by the negligibility of all addends in the decomposition above.


\end{document}